\documentclass[letterpaper, 10pt, journal]{IEEEtran}      % Use this line for a4
\usepackage{graphics} % for pdf, bitmapped graphics files
\usepackage{epsfig} % for postscript graphics files
\usepackage{epstopdf}
\usepackage{mathptmx} % assumes new font selection scheme installed
\usepackage{times} % assumes new font selection scheme installed
\usepackage{amsmath} % assumes amsmath package installed
\usepackage{amssymb}  % assumes amsmath package installed
\usepackage[noadjust]{cite}
\usepackage{filecontents}
\usepackage[usenames, dvipsnames]{color}
\usepackage{setspace}
\usepackage{psfrag}
\usepackage{bm}

\ifCLASSOPTIONcompsoc
  \usepackage[caption=false,font=normalsize,labelfont=sf,textfont=sf]{subfig}
\else
  \usepackage[caption=false,font=footnotesize]{subfig}
\fi

\usepackage[prependcaption,colorinlistoftodos,disable]{todonotes}

\newcommand{\todoiny}[1]{\todo[inline,color=yellow!60, linecolor=orange!250]{\bf\small#1}}

\newcommand{\il}[1]{{\color{black}{#1}}}
\newcommand{\vs}[1]{{\color{black}{#1}}}
\newcommand{\ill}[1]{{\color{black}{#1}}}
\newcommand{\illl}[1]{{\color{black}{#1}}}
\newcommand{\vsp}[1]{{\color{black}{#1}}}
\newcommand{\icl}[1]{{\color{black}{#1}}}

\newtheorem{theorem}{Theorem}
\newtheorem{remark}{Remark}
\newtheorem{assumption}{Assumption}
\newtheorem{definition}{Definition}
\newtheorem{lemma}{Lemma}
\newtheorem{prooft}{Proof of Theorem}
\newtheorem{proofl}{Proof of Lemma}
\usepackage{balance}

\title{\LARGE \bf
A system reference frame approach for \ill{stability}\\ analysis and control of power \ill{grids}}

\author{Chrysovalantis~Spanias and Ioannis~Lestas% <-this % stops a space
\thanks{C. Spanias is with the Distribution System Operator of the Electricity Authority of Cyprus and with the Department of Electrical Engineering, Computer Engineering and Informatics, Cyprus University of Technology. email: ca.spanias@edu.cut.ac.cy}% <-this % stops a space
\thanks{I. Lestas is with the \il{Department of Engineering, University of Cambridge, Trumpington Street, Cambridge CB21PZ, UK. email: icl20@cam.ac.uk}}}% <-this % stops a space
\begin{document}

\maketitle

%\thispagestyle{empty}
%\pagestyle{empty}

%%%%%%%%%%%%%%%%%%%%%%%%%%%%%%%%%%%%%%%%%%%%%%%%%%%%%%%%%%%%%%%%%%%%%%%%%%%%%%%%
\begin{abstract}
During the last decades, significant advances have been made in the area of power system stability and control.
%Although this development could assist the operation of the future power %grids, it
\ill{Nevertheless, when this analysis is carried out by means of decentralized conditions in a general network, it has been based}
%was based
on conservative assumptions such as the adoption of lossless networks. In the current paper, we present a novel %passivity-based
approach for decentralized stability analysis and control of power grids through the transformation of both the network and the bus dynamics into \ill{the} system reference frame. \ill{In particular}, the aforementioned transformation allows us to formulate the network model as \ill{an
%$(2 \times |\mathcal{N}|)$-input/$(2 \times |\mathcal{N}|)$
input-output system that is shown to be passive} even if the network's lossy nature is taken into account. We then introduce a broad class of bus dynamics
\ill{that are viewed as multivariable input/output systems compatible with the network formulation, and appropriate passivity conditions are imposed on those that guarantee stability of the power network.}
%that fits to the network's multi-variable formulation, and provide the %necessary passivity conditions that guarantee the asymptotic stability of the %whole power system in a decentralized manner.
We discuss \ill{the} opportunities \ill{and} advantages offered by this approach while explaining how \ill{this allows the inclusion of
advanced models for both generation and power flows.}
%various power system components, such as \ill{synchronous generators, and %dynamic loads.}
Our analysis is verified through applications \ill{to the} %generator buses of the
Two Area Kundur and the IEEE 68-bus test \ill{systems with both} primary frequency and voltage regulation mechanisms \ill{included}.
%\ill{via excitation control and power system stabilizers}.
\end{abstract}

\begin{IEEEkeywords}
power system stability, system reference-frame, passivity, frequency control, voltage control.
\end{IEEEkeywords}

%%%%%%%%%%%%%%%%%%%%%%%%%%%%%%%%%%%%%%%%%%%%%%%%%%%%%%%%%%%%%%%%%%%%%%%%%%%%%%%%
\section{INTRODUCTION} \label{sec:intro}

In the last decades, power systems have been through critical changes
%due to the undeniable need for mitigating the environmental pollution and
\ill{as a result of environmental issues and also to enhance} energy efficiency and security. Such changes are the introduction of new generation and storage technologies, and the rapid increase of the share of Renewable Energy Sources (RES) in power generation. Although, these advances contributed to \ill{technological and economic development}, they have \ill{resulted in a need} for increased stability, reliability and robustness \ill{in power} grids. Particularly, the large population of RES across the power networks, in combination with their intermittent nature affected \illl{the system} voltage and frequency stability \cite{kundur2004definition}.

This \illl{lack of effectiveness of existing} regulation mechanisms, coerced scientists to seek for new, more accurate and reliable control schemes. Additionally, several other \ill{highly distributed} power system components such as loads, inverter-based RES, Flexible Alternating Current Transmission System (FACTS) devices, \illl{were} also employed to provide fast-responding ancillary services to the grid and thus to assist in overcoming these voltage and frequency stability issues \vs{\cite{ma2013demand, kirby1997ancillary}}. However,
\ill{when a network wide stability analysis is carried out where one aims to deduce stability by means of local conditions, the complexity of the problem increases significantly and various simplifications are often employed in the literature that can hinder the applicability of the analysis. These include for example the assumption of a lossless network, independent study of voltage and frequency dynamics, or the lack of more advanced higher order dynamics in the feedback policies.}
%various
%the design and the justification of these control mechanisms were usually based on simple network models and bus dynamics in order to facilitate the stability analysis. Consequently, this led to the derivation of less accurate and more conservative control designs that do not always ensure power system stability and robustness \vs{\cite{stegink2016optimal}}.
%\todoiny{need some references in the paragraph above}

%Our approach in this paper is paper is to view a power %network

A key structural property that can facilitate power system stability \ill{analysis} %and help to  the aforementioned difficulties,
is the notion of passivity. %Passivity which has often been used in the analysis of %large-scale systems, can also allow the adoption of more %accurate, higher order dynamics and the deduction of %decentralized results when used appropriately %\cite{devane2016distributed}.
The application of passivity within power system studies dates back to the 80's, where passivity-based techniques were used to study the effect of Automatic Voltage Regulators (AVRs) on power systems \cite{miyagi1986stability}. More recently, the notion of passivity was widely used in power system studies via the framework of port-Hamiltonian systems (described in \cite{vanderschaft2013}). \ill{Examples of this approach include \cite{maschke2000energy, wang2003dissipative, fiaz2013port} as well as more recent works as in \cite{stegink2016unifying, schiffer2016stability, schiffer2014conditions, caliskan2014compositional, stegink2016optimal, andreasson2015distributed, trip2016internal, kasis2016primary}.
%. Three additional yet interesting, passivity-based %approaches can be found
What lies in common between the aforementioned studies is the fact that the stability analysis is carried out at each bus local machine reference frame. Specifically, both the bus and network dynamics are interfaced %within the elaborated analysis
through their expression at each bus dq-coordinates. Despite its %facile applicability and
broad use within the literature, this approach eliminates  natural passivity properties of the network and requires additional assumptions to be made to maintain those, such as that of a lossless network.}
%introduced several %issues %and
%disadvantages
%regarding the adoption of higher order dynamics, the incorporation of more %complex control mechanisms and the consideration of loss-less network models.}

In contrast to the related literature, this paper presents a
%broad %passivity-based
\ill{framework}
that facilitates the power system stability analysis through the transformation of both the network and the bus dynamics into \ill{the} system reference frame instead of each bus local reference frame. This transformation allows us to consider a lossy network model with arbitrary topology and to show that when \ill{an input/output} formulation is adopted, its passive nature is revealed. We then consider a broad class of bus dynamics \ill{that are viewed as multivariable input/output systems, compatible with the} %that fits to our multi-variable
network formulation and we provide \ill{appropriate local passivity conditions that ensure the asymptotic stability of the equilibrium points of the network}.
%of the interconnected system in a decentralized manner.
The aforementioned formulation \ill{provides the opportunity to incorporate into our analysis a variety of bus dynamics such as synchronous generators and loads,
and consider also frequency and voltage control policies. Throughout the paper we also show with realistic examples that the proposed conditions are satisfied by existing implementations when excitation control and \vsp{Power System Stabilizers (PSSs)} are present, and are hence not restrictive despite being decentralized.}
%inverter-based RES, FACTS devices etc. and to adopt higher order dynamic models %that could enhance the reliability and the robustness of the power system.

The \ill{paper} is organized as follows:
In Section \ref{sec:background}, we introduce some basic preliminaries regarding power system analysis \ill{that will be used within the paper}. The formulation of the proposed approach is presented in Section \ref{sec:formulation}. In Section \ref{sec:discuss}, we provide an extensive discussion \ill{on the} advantages and \illl{opportunities} provided by the proposed framework. Section \ref{sec:simulations} illustrates our results through simulations on the Kundur 2-area and the IEEE 68-bus test systems. Finally, conclusions are drawn in Section~\ref{sec:conclusions}.

\vspace{-1mm}
%%%%%%%%%%%%%%%%%%%%%%%%%%%%%%%%%%%%%%%%%%%%%%%%%%%%%%%%%%%%%%%%%%%%%%%%%%%%%%%%
\section{BACKGROUND} \label{sec:background}

\ill{The analysis framework that will be presented in this paper relies on the representation of power systems as an interconnection of dynamical systems in an appropriate frame of reference. In order to describe this framework we present in this section some preliminaries that are relevant in this context.}

%%%%
\vspace{-2mm}
\subsection{Alternating Current (AC) three-phase sources}
%Power systems are networks consisting of devices that generate, transmit and distribute electrical %energy to consumers. The majority of power systems and their components rely upon three-phase AC power \cite{ arthur2000power}.
 \ill{We will use the notation  %will be represented in the paper as}
\vspace{-1mm}
\begin{equation*}
x_{ABC}= [ x_{A}(t) \ x_{B}(t) \ x_{C}(t) ]^\textrm{T}
\end{equation*}
to represent three-phase AC signals \illl{$x_{ABC}:\mathbb{R}^+ \to \mathbb{R}^3$}.}
%where $x_{ABC}:\mathbb{R}^+ \to \mathbb{R}^3$ is a three-phase AC signal.
\ill{In particular, three-phase voltages and currents will be denoted~as}
%Thus, the three-phase voltages and currents can be written as
\begin{equation} \label{eq:voltcurr}
v_{ABC}= [ v_{A}(t) \ v_{B}(t) \ v_{C}(t) ]^\textrm{T} \ \mbox{and} \ i_{ABC}= [ i_{A}(t) \ i_{B}(t) \ i_{C}(t) ]^\textrm{T}
\end{equation}
respectively.

\begin{assumption}\label{assumption:balanced}
\ill{The power networks that will be considered in the paper} consist of symmetric, balanced, positive-sequence, three-phase AC generation sources.
\end{assumption}

Since power systems are designed to be symmetric and balanced, the above assumption is \ill{often} accurate, especially when analysis is carried out at \ill{the} transmission level.  Assumption \ref{assumption:balanced} results in three symmetric waveforms which have $120^{o}$ phase difference between each other, i.e.

\begin{equation}  \label{eq:waves}
x_{ABC}=
\begin{bmatrix}
x_{A}(t) \\
x_{B}(t) \\
x_{C}(t)
\end{bmatrix}
=\sqrt{2}|x|
\begin{bmatrix}
\cos(\gamma_x (t)) \\
\cos(\gamma_x (t)-\frac{2\pi}{3}) \\
\cos(\gamma_x (t)+\frac{2\pi}{3})
\end{bmatrix}
\end{equation}
where $|x|\in \mathbb{R}^+$ is \ill{the} amplitude and $\gamma_x \in [0,2\pi)$ is \ill{the} phase of the waveform. The fact that the three phases are balanced results in
\begin{equation} \label{eq:neutral}
x_{A}(t) + x_{B}(t) + x_{C}(t) = 0
\end{equation}
\ill{Furthermore}, problems in symmetric and balanced power systems can be dealt \ill{with} by using only the phase A and then \ill{deduce the results} for phases B and C \ill{from} (\ref{eq:waves}).

%%%%
\subsection{Phasor representation}

To simplify power system analysis, it is usually convenient to use the phasor representation of voltages and currents rather than their sinusoidal form (\ref{eq:waves}). The phasor representation is defined as follows \vs{\cite{glover2012power}}:

\begin{definition} \label{definition:phasors}
A phasor is a complex number representing a sinusoidal signal
\begin{equation} \label{eq:sinusoid}
x(t)=|x|\cos(\gamma_x (t))=|x|\cos(\omega t + \phi_x)
\end{equation}
whose amplitude $|x|$, frequency $\omega$ and phase angle $\phi_x$ can be time varying quantities. Using the quantity $\bar{X}$ to indicate the phasor, the polar phasor representation of the signal (\ref{eq:sinusoid}) is given by:
\begin{equation} \label{eq:polar}
\bar{X}=|x| e^{\textrm{j} \gamma_x (t)} =|x| \angle{\gamma_x (t)}.
\end{equation}
We can also obtain its rectangular representation by using Euler's identity as follows:
\begin{equation} \label{eq:rectangular}
\bar{X}=|x| e^{\textrm{j} \gamma_x (t)} =|x| \big(\cos( \gamma_x (t)) + \textrm{j} \sin( \gamma_x (t)) \big).
\end{equation}
\vs{\illl{A representation often adopted is to have a constant}
$\omega = \omega_s = 2 \pi f_s$ where $f_s$ denotes the synchronous frequency of a power grid (50 or 60 Hz), \illl{and represent phasors as:} %the phasor representation (\ref{eq:polar}) \illl{can be reduced} to:
\begin{equation} \label{eq:phasor}
\bar{X}=|x| e^{\textrm{j} \phi_x} =|x| \angle{\phi_x}.
\end{equation}}
\hspace{-0.2cm}Note that $\phi_x$ can be a time varying quantity that models variations in frequency.
\end{definition}
%The majority of the related literature (\cite{kundur1994power, machowski2011power, arthur2000power, sauer1997power}), defines a phasor as a complex sinusoidal quantity with constant frequency $\omega = \omega_s = 2 \pi f_s$ where $f_s$ denotes the synchronous frequency of the network (50 or 60 Hz).
%\todoiny{Is this true, e.g. does the bergen and vittal book use only the synchronous frequency?}
%The phasor representation (\ref{eq:polar}) therefore reduces to:
%\begin{equation} \label{eq:phasor}
%\bar{X}=|x| e^{j \phi_x} =|x| \angle{\phi_x}.
%\end{equation}
%Similarly to \cite{anderson2008power}, in this paper we consider that frequency $\omega$ varies with time and we adopt the representation (\ref{eq:polar}).
%%%%
\vspace{-.4cm}
\subsection{(0,d,q) or Park's transformation}

A key tool to facilitate power systems analysis is $(0,d,q)$ or Park's transformation. The sinusoidal waveforms (\ref{eq:waves}), describing either voltages or currents,  introduce significant complexity in
\ill{the analysis.
%solving several power system problems due to their time dependency.
Therefore,} to simplify these equations, we use \illl{the} $(0,d,q)$ or Park's transformation so as to map \illl{the} system's components into three axis that rotate at a specific velocity $\omega$, namely, the $0$-axis, the $d$-axis and the $q$-axis. Following \cite{kundur1994power, machowski2011power, arthur2000power, sauer1997power}, the Park's transformation is defined by:
\begin{equation} \label{eq:transf1}
\small{
\begin{bmatrix}
x_0 \\
x_d \\
x_q
\end{bmatrix}
= \sqrt{\frac{2}{3}}
\underbrace{\begin{bmatrix}
\frac{1}{\sqrt{2}} & \frac{1}{\sqrt{2}} & \frac{1}{\sqrt{2}} \\
\cos \rho{(t)} & \cos{(\rho (t)-\frac{2\Pi}{3})} & \cos{(\rho (t)+\frac{2\Pi}{3})} \\
\sin{\rho (t)} & \sin{(\rho (t)-\frac{2\Pi}{3})} & \sin{(\rho (t)+\frac{2\Pi}{3})}
\end{bmatrix}}_{\displaystyle P}
\begin{bmatrix}
x_A \\
x_B \\
x_C
\end{bmatrix}}
\end{equation}
\il{where $P$ is the transformation matrix relating the $abc$ and $0dq$ vectors.}
%which defines the $P$ transformation matrix and the $abc$ and $0dq$ vectors.
The new $0dq$ variables are also called Park's variables. Furthermore, the Park's transformation is orthogonal, i.e. $P^{-1}=P^\textrm{T}$. Under Assumption \ref{assumption:balanced}, which \il{yields} the zero sum of both the voltages and \illl{currents} of the three phases (equation (\ref{eq:neutral})), the $0$-component in (\ref{eq:transf1}) is equal to zero and can be therefore neglected. Now, considering that the $0$-component \illl{can be} neglected, we substitute equation (\ref{eq:waves}) into (\ref{eq:transf1}) to get
\vsp{\begin{equation} \label{eq:transf4}
x_{dq}=\sqrt{3} |x|
\begin{bmatrix}
\cos{(\gamma_x (t) - \rho (t))} \\
\sin{(\gamma_x (t) - \rho (t))}
\end{bmatrix}
\end{equation}}
which is essentially a projection of phasors onto axes rotating with frequency $\omega=\dot{\rho}$.  Similarly to the \illl{$abc$} components, the $dq$ components can be also expressed as complex numbers onto these rotating axes, i.e. $X_{dq}=X_q + \textrm{j} X_d$. \illl{This representation will be referred to as the phasor representation of $x$ in a frame of reference rotating with frequency $\dot{\rho}$.

Note also that $\bar X$ in \eqref{eq:phasor} is the phasor representation in a frame of reference rotating with a constant frequency $\omega_s$. The latter will be referred to as the system reference frame.}
%%%%
\subsection{Power Network}

A power network with arbitrary topology can be described by a connected and undirected graph $(\mathcal{N},\mathcal{E})$, where $\mathcal{N}=\{ 1, 2, \ldots |\mathcal{N}|\}$ is the set of buses and $\mathcal{E} \subset \mathcal{N} \times \mathcal{N}$ the set of transmission lines connecting the buses.  We use $(i,j)$ to denote the link connecting the network buses $i$ and $j$. Based on the formulation described in \cite{schiffer2016survey}, we consider the following assumptions in order to derive the equations describing the network.

\begin{assumption} \label{assumption:rlc}
Transmission lines can be represented by symmetric three-phase RLC elements.
\end{assumption}

\begin{assumption} \label{assumption:dynamics}
Transmission line dynamics evolve on a much faster timescale than the dynamics of the generation sources and the loads.
\end{assumption}

Assumption \ref{assumption:dynamics} states that transmission lines reach steady state much earlier than the generators and the loads. Thus, the power network can be modeled by the static network current flows given by the nodal set of equations:
\begin{equation}\label{eq:network}
\bar{I}=Y_n \bar{V} = ( G_n + \textrm{j} B_n ) \bar{V}.
\end{equation}
$Y_n \in\mathbb{C}^{|\mathcal{N}| \times |\mathcal{N}|}$, and $G_n$, $B_n \in \mathbb{R}^{|\mathcal{N}| \times |\mathcal{N}|}$ are the network's admittance, conductance and susceptance matrices respectively. $\bar{I} \in \mathbb{C}^{|\mathcal{N}|}$ and  $\bar{V} \in \mathbb{C}^{|\mathcal{N}|}$ denote the net injected current and the bus voltage vectors of the power grid \ill{in their phasor representation}. The derivation of the nodal admittance matrix (\ref{eq:network}) is extensively described in \cite{arthur2000power} and is based on the fact that the transmission lines are modeled by their $\Pi$-equivalent \il{model} according to Assumption \ref{assumption:rlc}.

\begin{remark}
$G_n$ and $B_n$ are real, $|\mathcal{N}| \times |\mathcal{N}|$, sparse symmetric matrices and they do not include \ill{loads or} FACTS devices and line compensation components.
\end{remark}

The components of net injected current and bus voltage vectors can equivalently be expressed in \ill{either their rectangular or} polar complex form. However,
%in order to develop the analytic network equations,
\ill{it is convenient here to express these} in the same form as the elements of the nodal admittance matrix, that is, the rectangular form. Considering \ill{a steady} state network frequency \ill{$\omega_s$}, \ill{the} net injected currents and bus voltages can be written \illl{using the phasor representation in \eqref{eq:phasor} as}
\begin{equation}\label{eq:current}
\bar{I}_i = I_i \angle \phi_{I,i} = I_i \cos{\phi_{I,i}} + \textrm{j} I_i\sin{\phi_{I,i}} = I_{a,i} + \textrm{j} I_{b,i}
\end{equation}
% and
\begin{equation}\label{eq:voltage}
\bar{V}_i = V_i \angle \phi_{V,i} = V_i \cos{\phi_{V,i}} + \textrm{j} V_i\sin{\phi_{V,i}} = V_{a,i} + \textrm{j} V_{b,i}
\end{equation}
respectively, for all $i \in \mathcal{N}$. We now define the vectors $I_a = [I_{a,1} \ I_{a,2} \ ... \ I_{a,|\mathcal{N}|}]^\textrm{T}$, $I_b = [I_{b,1} \ I_{b,2} \ ... \ I_{b,|\mathcal{N}|}]^\textrm{T}$, $V_a = [V_{a,1} \ V_{a,2} \ ... \ V_{a,|\mathcal{N}|}]^\textrm{T}$ and $V_b = [V_{b,1} \ V_{b,2} \ ... \ V_{b,|\mathcal{N}|}]^\textrm{T}$ $\in \mathbb{R}^{|\mathcal{N}|}$. The net injected current and the bus voltage vectors can therefore be written as:
\begin{equation}\label{eq:currentvoltagevec}
\bar{I}= I_a + \textrm{j} I_b \ \mbox{and} \ \bar{V}= V_a + \textrm{j} V_b
\end{equation}
respectively. By substituting equations (\ref{eq:currentvoltagevec}) into (\ref{eq:network}) we get:
\begin{equation}\label{eq:network2}
\bar{I} = I_a + \textrm{j} I_b =(G_n V_a - B_n V_b) + \textrm{j} (B_n V_a + G_n V_b).
\end{equation}
\ill{Equation (\ref{eq:network2}) is then used to deduce the equations} for the net injected current components, $I_{a,i}$  and $I_{b,i}$, at each bus $i = 1,2, \dots , |\mathcal{N}|$, \ill{i.e.,} we get:
\begin{equation}\label{eq:networkcomps}
\small{
I_{a,i} = \sum_{j=1}^{|\mathcal{N}|} (G_{ij} V_{a,j} - B_{ij} V_{b,j}) \mbox{ and } I_{b,i} = \sum_{j=1}^{|\mathcal{N}|} (B_{ij} V_{a,j} + G_{ij} V_{b,j})}.
\end{equation}

\illl{Note that in the network equations \vs{\eqref{eq:network} - \eqref{eq:networkcomps}},  the current and voltage phasors $\bar I, \bar V$ are represented in the \textit{system reference frame}, i.e. a common reference frame rotating at the synchronous frequency $\omega_s$. The network admittance matrix in \eqref{eq:network} is also evaluated at $\omega_s$. This is a common approach in the literature, and as discussed in \cite{schiffer2016survey}, it is a valid approximation, %even when the frequency is time varying,
under the assumption that transmission line dynamics are much faster than machine~dynamics.

%\illl{When the network frequency is time varying a common approach in the literature is to consider a system reference frame for the network equations. In particular, it is assumed that the network equations  \eqref{eq:networkcomps}, \eqref{eq:network} are valid when the current and voltage phasors $\bar I, \bar V$ in these equations  denote phasors in the %the $(a,b)$ or the
%\textit{system reference frame}, i.e. a common reference frame rotating at the synchronous frequency $\omega_s$. The network admittance matrix in \eqref{eq:network} is also evaluated at $\omega_s$. As discussed in \cite{schiffer2016survey},  this is a valid approximation, commonly used in the literature, under the assumption that transmission line dynamics are much faster than machine~dynamics\footnote{\ill{In the remainder of the paper we will use the notation  $\bar{I}= I_a + j I_b, \bar{V}= V_a + j V_b$ to denote phasor representations in the system reference frame.}}.
%The  network equations \ill{(\ref{eq:networkcomps}), \eqref{eq:network}} are \ill{usually expressed in the literature} in the $(a,b)$ or the \textit{system reference frame}, \ill{i.e. a common reference frame rotating at the synchronous frequency $\omega_s$. The network admittance matrix in \eqref{eq:network} is also evaluated at $\omega_s$. As discussed in \cite{schiffer2016survey},  the latter is a valid approximation, commonly used in the literature, also when the frequency of the signals is time varying, under the assumption that transmission line dynamics are much faster than machine~dynamics.

It is} also important to \ill{consider} the transition from the system reference frame to the \ill{local $(d,q)$}  or the \textit{machine reference frame} and vice versa. We thus define the angle $\delta_i \in [0,2\pi)$ denoting the phase difference between the local machine reference frame at bus $i$, with phase angle $\rho_i (t)$, and the system reference frame which rotates at synchronous frequency $\omega_s$, i.e.
\begin{equation}  \label{eq:angle}
\delta_i = \int_0^t \il{(\dot{\rho}_i (\tau) - \omega_s)} \ d\tau \ \Rightarrow \ \dot{\delta}_i = \dot{\rho}_i (t) - \omega_s = \omega_i - \omega_s
\end{equation} %We note here that in synchronous generator dynamics, the phase angle $\rho_i (t)$ denotes the phase angle of the stator's voltage $\bar{E}_i$.

\begin{figure}[!t]
\centering
\subfloat[The geometry of $dq$-components on the system's complex plane.]{\includegraphics[width=0.47\columnwidth]{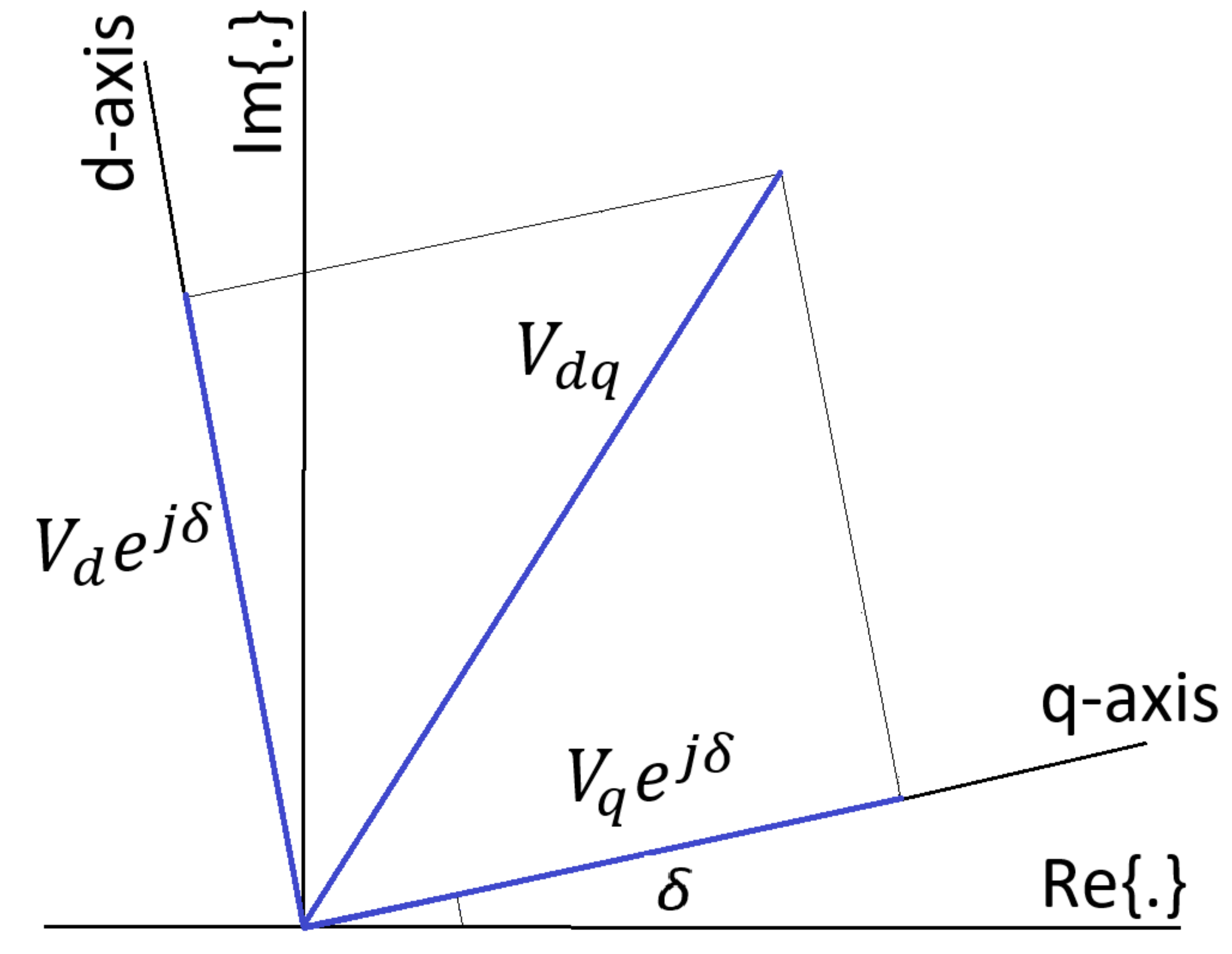}
\label{fig:referenceframes1}}
\hfil
\subfloat[The decomposition of $dq$-components in both the system and the machine reference frames.]{\includegraphics[width=0.47\columnwidth]{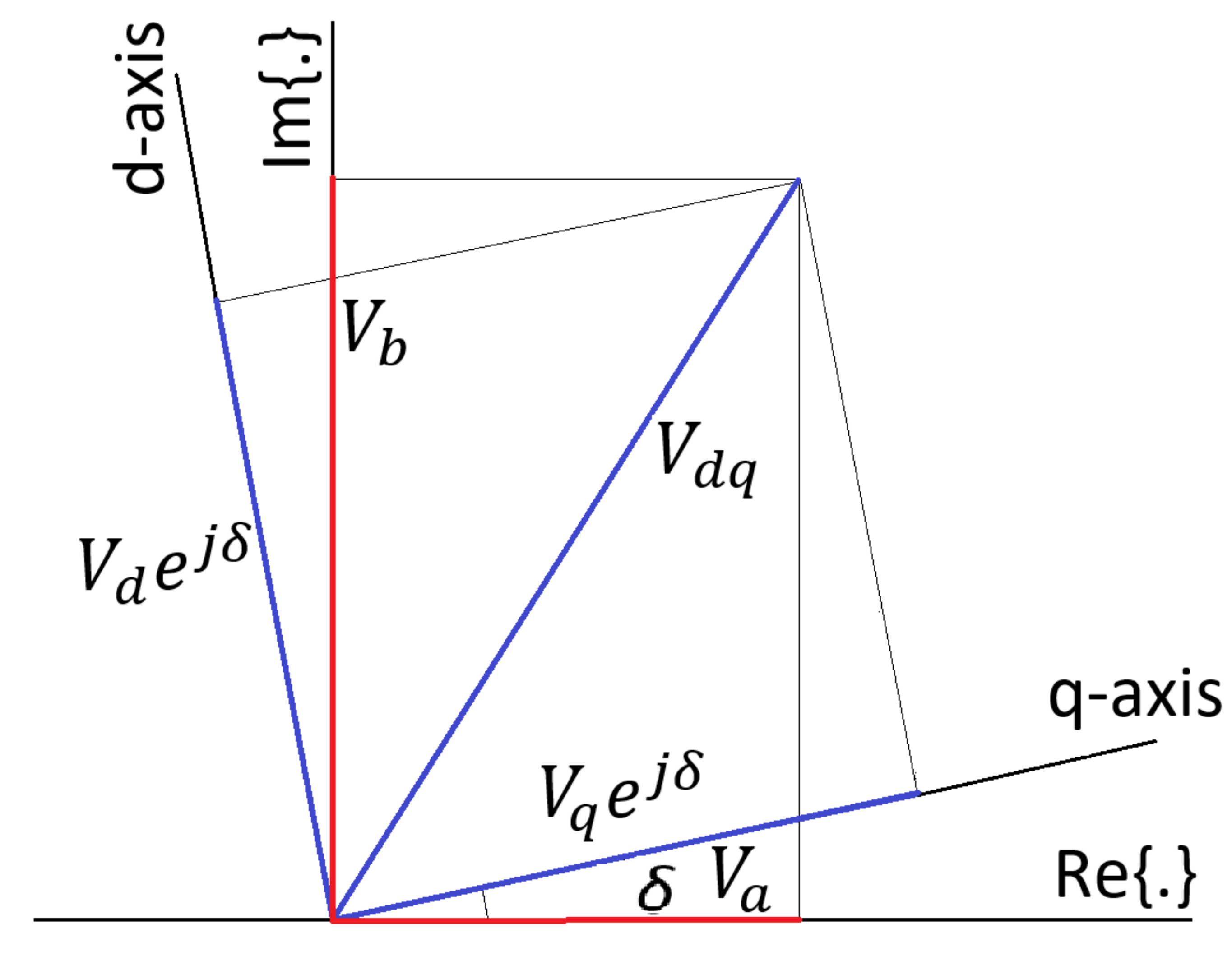}
\label{fig:referenceframes2}}
\caption{Relative position of the machine reference frame with respect to the system reference frame \cite{arthur2000power}.}
\label{fig:referenceframes}
\end{figure}

The relative position of the two systems of coordinates is illustrated in Figure \ref{fig:referenceframes} and the relationship between them is given by:
\begin{equation}  \label{eq:mapping}
V_{dq,i}=T(\delta_i) \bar{V}_i \Leftrightarrow
\begin{bmatrix}
V_{q,i} \\
V_{d,i}
\end{bmatrix}
=
\begin{bmatrix}
\cos{\delta_i} & \sin{\delta_i} \\
-\sin{\delta_i} & \cos{\delta_i}
\end{bmatrix}
\begin{bmatrix}
V_{a,i} \\
V_{b,i}
\end{bmatrix}
\end{equation}
where the transformation matrix $T(\delta_i)$ denotes \vs{the mapping of the phasor components in the system reference frame to the $dq$-components for bus $i$}. \todoiny{isn't this the other way round, i.e. $T(\delta_i)$ maps the phasor components in system reference frame to the dq components}The transformation $T$ is also orthogonal ($T^{-1}=T^\textrm{T}$), and its inverse transformation can be written as:
\begin{equation}  \label{eq:invmapping}
\bar{V}_i =T^{-1}(\delta_i) V_{dq,i} \Leftrightarrow
\begin{bmatrix}
V_{a,i} \\
V_{b,i}
\end{bmatrix}
=
\begin{bmatrix}
\cos{\delta_i} & -\sin{\delta_i} \\
\sin{\delta_i} & \cos{\delta_i}
\end{bmatrix}
\begin{bmatrix}
V_{q,i} \\
V_{d,i}
\end{bmatrix}.
\end{equation}

Equivalently, for the net current injection components we get
\begin{equation}  \label{eq:mapping2}
I_{dq,i}=T(\delta_i) \bar{I}_i
%\Leftrightarrow
%\begin{bmatrix}
%I_{q,i} \\
%I_{d,i}
%\end{bmatrix}
%=
%\begin{bmatrix}
%\cos{\delta_i} & \sin{\delta_i} \\
%-\sin{\delta_i} & \cos{\delta_i}
%\end{bmatrix}
%\begin{bmatrix}
%I_{a,i} \\
%I_{b,i}
%\end{bmatrix}
\end{equation}
\begin{equation}  \label{eq:invmapping2}
\bar{I}_i =T^{-1}(\delta_i) I_{dq,i}
%\Leftrightarrow
%\begin{bmatrix}
%I_{a,i} \\
%I_{b,i}
%\end{bmatrix}
%=
%\begin{bmatrix}
%\cos{\delta_i} & -\sin{\delta_i} \\
%\sin{\delta_i} & \cos{\delta_i}
%\end{bmatrix}
%\begin{bmatrix}
%I_{q,i} \\
%I_{d,i}
%\end{bmatrix}
\end{equation}

%\begin{figure}[!t]
%\centering
%\begin{subfigure}{0.24\textwidth}
%\centering
%\includegraphics[width=1.0\linewidth]{ref1}
%\caption{The geometry of $dq$-components on the system's complex plane.}
%\label{fig:referenceframes1}
%\end{subfigure}
%\begin{subfigure}{0.24\textwidth}
%\centering
%\includegraphics[width=1.0\linewidth]{ref2}
%\caption{The decomposition of $dq$-components in both the system and the machine reference frames.}
%\label{fig:referenceframes2}
%\end{subfigure}
%\caption{Relative position of the machine reference frame with respect to the system reference frame \cite{arthur2000power}.}
%\label{fig:referenceframes}
%\end{figure}

Now, we express network equations (\ref{eq:network2}) in each generator's reference frame so as to obtain the general network relationships. By substituting (\ref{eq:invmapping}) and (\ref{eq:invmapping2}) into (\ref{eq:networkcomps}) we get:
\begin{equation}\label{eq:networkcompsdq}
\begin{split}
I_{q,i} =& \sum_{j=1}^{|\mathcal{N}|} \Big[ V_{q,j} \Big( G_{ij} \cos(\eta_{ij}) + B_{ij} \sin(\eta_{ij})\Big) \\ &+ V_{d,j} \Big( G_{ij} \sin(\eta_{ij}) - B_{ij} \cos(\eta_{ij}) \Big) \Big] \\
I_{d,i} =& \sum_{j=1}^{|\mathcal{N}|} \Big[ V_{q,j} \Big( -G_{ij} \sin(\eta_{ij}) + B_{ij} \cos(\eta_{ij}) \Big) \\&+ V_{d,j} \Big( G_{ij} \cos(\eta_{ij}) + B_{ij} \sin(\eta_{ij}) \Big) \Big]
\end{split}
\end{equation}
where, for ease of notation, angle differences are written as $\eta_{ij} = \delta_i - \delta_j$.

\section{FRAMEWORK FORMULATION} \label{sec:formulation}

%In the current section we present our framework formulation. We first form the network model as a $(2 \times |\mathcal{N}|)$-input/$(2 \times |\mathcal{N}|)$-output system, and show that under such formulation the power networks are passive systems. We then introduce a broad class of bus dynamics which, in contrast to the recent literature, are now transformed into system reference frame. The transformation of bus dynamics from their loqal dq-reference frame into the system reference frame is carried out by incorporating the mappings (\ref{eq:mapping})-(\ref{eq:invmapping2}) into their dynamic equations. The following subsections deal with the quantification of what is meant by an equilibrium of the interconnected system, and the presentation of the necessary passivity conditions that have to be satisfied by bus dynamics so as to guarantee the stability of the aforementioned equilibria.

%%%%
\subsection{Network Model}
\label{sec:NetworkModel}
\begin{figure}[t!]%[htbp!]
\centering
\includegraphics[width=\columnwidth]{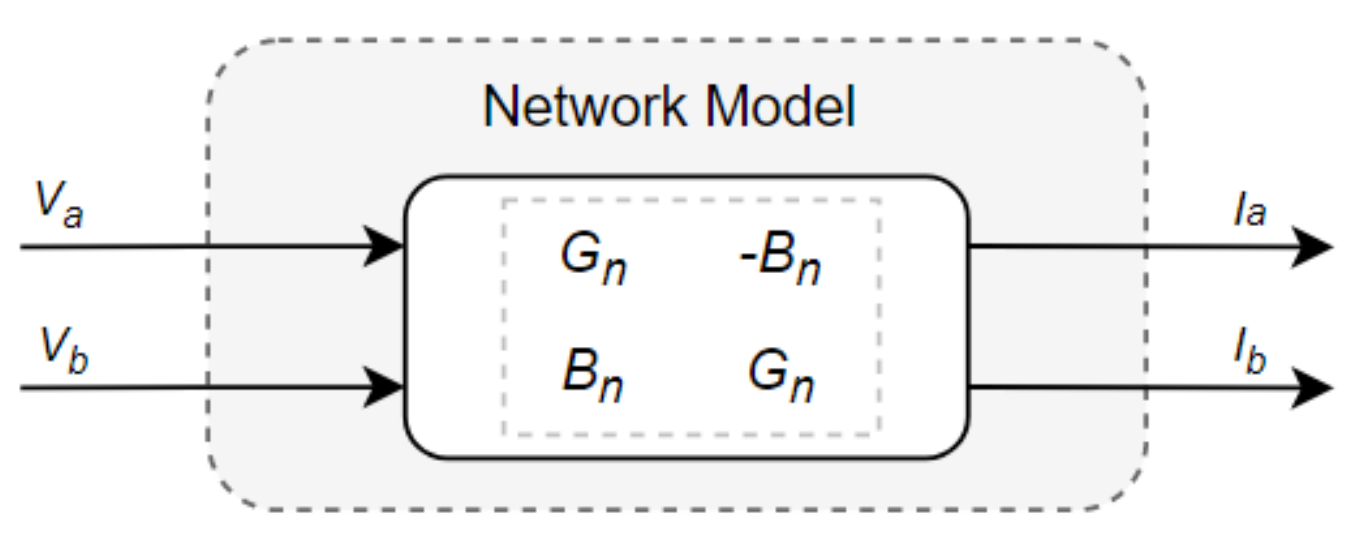}
\caption{The $(2 \times |\mathcal{N}|)$-input/$(2 \times |\mathcal{N}|)$-output system that is used to model the power network.}
\label{fig:networkmodel}
\end{figure}

As discussed in Section \ref{sec:background}, the power network is represented by the nodal set of equations (\ref{eq:network}), which can be also written in the \ill{rectangular} form (\ref{eq:network2}). In order to formulate the network model, we separate the real and the imaginary part \ill{of equation} (\ref{eq:network2}) so as to form the following $(2 \times |\mathcal{N}|)$-input/$(2 \times |\mathcal{N}|)$-output system
\begin{equation} \label{eq:network4}
\begin{bmatrix} I_a \\ I_b \end{bmatrix} = \begin{bmatrix} G_n  & -B_n \\ B_n & G_n \end{bmatrix} \begin{bmatrix} V_a \\ V_b \end{bmatrix} \\  = H_{2n} \begin{bmatrix} V_a \\ V_b \end{bmatrix} = g^N ([V_a^\textrm{T} \ V_b^\textrm{T}])
\end{equation}
where $H_{2n}$ denotes the matrix relating the vectors $[ V_a^\textrm{T} \ V_b^\textrm{T} ]^\textrm{T}$ with the vectors $[ I_a^\textrm{T} \ I_b^\textrm{T} ]^\textrm{T}$. The vector function \ill{$g^N: \ \mathbb{R}^{2|\mathcal{N}|} \to \mathbb{R}^{2|\mathcal{N}|}$} provides an alternative \ill{notation}
% to represent} the network system
so as to comply with the definitions that we are about to use in the forthcoming paragraphs. The aforementioned system is illustrated in Figure~\ref{fig:networkmodel}.

We are now ready to examine the passivity properties that are revealed through the aforementioned modeling. Taking into account \il{that} Assumption \ref{assumption:dynamics} holds, we first provide the following fundamental passivity definition \cite{Kha01}.
\begin{definition} \label{definition:passivity1}
Consider the system described by the memory-less function $y=g(t,u)$ where $g:[0,\infty) \times \mathbb{R}^p \to \mathbb{R}^p$. This system is passive if $u^\textrm{T} y \geq 0$.
\end{definition}

As stated above, the static network model (\ref{eq:network4}) is passive if and only if the inequality within Definition \ref{definition:passivity1} is satisfied, i.e.
\begin{equation} \label{eq:condition1}
u^\textrm{T} y = [V_a^\textrm{T} \ V_b^\textrm{T}] \begin{bmatrix} I_a \\ I_b \end{bmatrix} \geq 0
\end{equation}
for all $V_a$, $V_b$, $I_a$, $I_b$ $\in \illl{\mathbb{R}^{|\mathcal{N}|}}$.

\begin{lemma} \label{lemma:passivity}
The network system defined in (\ref{eq:network4}) with inputs the vectors of bus voltage components $[ V_a^\textrm{T} \ V_b^\textrm{T} ]^\textrm{T}$ and outputs the vectors of net injected current components $[ I_a^\textrm{T} \ I_b^\textrm{T} ]^\textrm{T}$ is passive.
\end{lemma}

\begin{remark}
We see within the proof of Lemma \ref{lemma:passivity} \ill{(provided in the appendix)} that \il{condition} (\ref{eq:condition1}) always holds and that the passivity of the network system is ensured regardless of its topology. Specifically, due to the form of the composite matrix $H_{2n}$, the positive semi-definiteness of the network's conductance matrix $G_n$ is sufficient \illl{for condition} (\ref{eq:condition1}) to be satisfied. $G_n$ in turn, is always positive semi-definite since it has positive diagonal elements and \il{is diagonally dominant.}
\end{remark}

\begin{remark}
As we are about to discuss in the next section, the majority of the recent literature dealing with power system stability \ill{in general network topologies} adopts lossless networks, i.e. $G_n = 0$. The main reason for considering such simplification lies in the fact that when the analysis is carried out in $dq$-coordinates, the passivity \il{property} holds only for lossless networks. For the proposed approach, under such assumption \illl{condition} (\ref{eq:condition1}) becomes
\begin{equation} \label{eq:condition3}
u^\textrm{T} y = [V_a^\textrm{T} \ V_b^\textrm{T}]  \begin{bmatrix} 0 & -B_n \\ B_n & 0 \end{bmatrix} \begin{bmatrix} V_a \\ V_b \end{bmatrix} = 0
\end{equation}
\illl{Note} that the network's passivity \illl{follows here easily from the skew-symmetry of the matrix~$H_{2n}$}.
\end{remark}

%%%%
\subsection{Bus Dynamics}

\begin{figure}[t!]%[htbp!]
\centering
\includegraphics[width=\columnwidth]{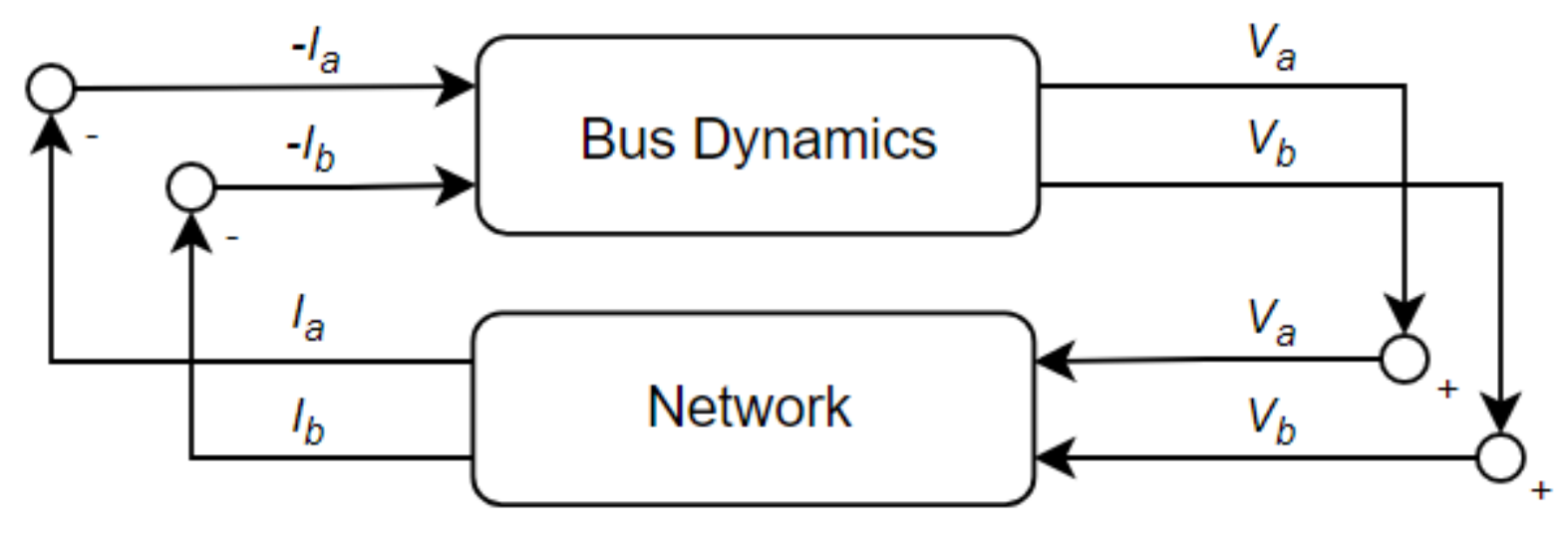}
\caption{\il{The power network represented as an interconnection of input/output systems associated with the bus dynamics and transmission lines, respectively.}}
%network system interconnected as a negative feedback to the network buses.}
\label{fig:networkA}
\end{figure}

In order to \ill{incorporate} %\ill{connect} the network with
the bus models and derive %several
\il{stability} results for the interconnected system, both the network and the bus dynamics have to be expressed in the same reference frame, \illl{which is chosen here as} the {system reference frame}. In contrast to the recent literature, we therefore transform the bus dynamics into \ill{the} system reference frame instead of each bus local $dq$-coordinates, and consider that each of the $|\mathcal{N}|$ buses forms a 2-input/2-output system so as to fit with the \ill{network formulation described in the previous section}. A graphical representation of the interconnected system is provided in Figure \ref{fig:networkA} where the multi-input/multi-output network system is \il{connected} %as a negative feedback
to the aggregate bus \il{dynamics}. The bus models are expressed in \illl{the} system reference frame by incorporating the mappings $T$ and $T^{-1}$ (equations (\ref{eq:mapping})-(\ref{eq:invmapping2})) into \ill{the} bus dynamics. This approach is to the best of our knowledge novel,  \il{and} allows the consideration of more relaxed conditions for the network while giving the opportunity for decentralized stability analysis and control.

%By following the formulation presented in the two-part paper %\cite{kasis2016primary, devane2016primary},
We now introduce a broad class of systems that are used to represent the bus models. We consider that these dynamical systems have inputs the phasor components of the net current injection $(-I_{a,i}, -I_{b,i}) \in \mathbb{R}^2$, states $x_i \in X \subseteq \mathbb{R}^k$ and outputs the phasor components of the bus voltage $(V_{a,i}, V_{b,i}) \in \mathbb{R}^2$. The state-space representation of the aforementioned dynamical systems is given by:\il{
\begin{equation} \label{eq:model}
\begin{split}
\dot{x}_i =& f_i (x_i, u_i) \\%\quad i \in \mathcal{N} \\
y_i %=& \begin{bmatrix} V_{a,i} \\ V_{b,i} \end{bmatrix}
=& g_i (x_i, u_i) \quad i \in \mathcal{N}
\end{split}
\end{equation}
where $u_i=[-I_{a,i}, -I_{b,i}]$, $y_i=[V_{a,i}, V_{b,i}]$.}
The vector functions \ill{$f_i : \mathbb{R}^{k_i} \times \mathbb{R}^{2} \to \mathbb{R}^{k_i}_i$ and $g_i : \mathbb{R}^{k_i} \times \mathbb{R}^{2} \to \mathbb{R}^2$} are locally Lipschitz for any \ill{$i\in \mathbb{N}$}. We also note here that the bus dynamics (\ref{eq:model}) can be of arbitrary dimension.

We \ill{now} describe \ill{what is meant by an equilibrium} of the interconnected system (\ref{eq:network4}) \illl{and (\ref{eq:model}).}

\begin{definition} \label{definition:equilibria}
\il{The constant vector
%} $(\hat{I}_a , \hat{I}_b, \hat{x}, \hat{V}_a, \hat{V}_b)$ define
$\hat x=[\hat x_1 \ \hat x_2 \ \hdots \ \hat{x}_{|\mathcal{N}|}]$, $\hat x_i\in\mathbb{R}^{k_i}$ \ill{is} an} equilibrium of the interconnected system (\ref{eq:network4}) and (\ref{eq:model}), if the time derivative of the states \illl{in (\ref{eq:model})} is equal to zero when\footnote{\illl{Note that the inputs $u_i$ are also function of the states $x_1, \hdots, x_{|\mathcal{N}|}$ as the system is interconnected.}} \illl{$x_i=\hat x_i, \  i \in \mathcal{N}.$}
%\begin{equation} \label{eq:equilibrium}
%$\ill{f_i(\hat x_i)} = 0 \ \  i \in \mathcal{N}.$}
%\end{equation}
\end{definition}

%Since we have introduced both the network and the bus models, we can now define the equilibria of the interconnected system (\ref{eq:network4}) and (\ref{eq:aggr_model}).
%\begin{definition} \label{definition:equilibria}
%Each of the constants $(\hat{I}_a , \hat{I}_b, \hat{x}, \hat{V}_a, \hat{V}_b)$ define an equilibrium of the interconnected system (\ref{eq:network4}) and (\ref{eq:aggr_model}), if the following hold
%\begin{equation} \label{eq:equilibrium}
%\begin{split}
%0 =& f_i (\hat{x}_i, -\hat{I}_{a,i}, -\hat{I}_{b,i}) \quad i \in \mathcal{N}\\
%\begin{bmatrix} \hat{V}_{a,i} \\ \hat{V}_{b,i} \end{bmatrix} =& g_i (\hat{x}_i , -\hat{I}_{a,i}, -\hat{I}_{b,i}) \quad i \in \mathcal{N}\\
%\hat{I}_{a,i} =& \sum_{j=1}^{|\mathcal{N}|} (G_{ij} \hat{V}_{a,j} - B_{ij} \hat{V}_{b,j}) \quad (i,j) \in \mathcal{E} \\
%\hat{I}_{b,i} =& \sum_{j=1}^{|\mathcal{N}|} (B_{ij} \hat{V}_{a,j} + G_{ij} \hat{V}_{b,j}) \quad (i,j) \in \mathcal{E}
%\end{split}
%\end{equation}
%We call (\ref{eq:equilibrium}) the equilibrium conditions for the closed loop system (\ref{eq:network4}) and (\ref{eq:aggr_model}).
%\end{definition}

%%%%
\subsection{Passivity Conditions on Bus Dynamics}

%Based on \cite{kottenstette2010relationships},
\ill{We now present passivity conditions %which when satisfied by
on the bus dynamics, which when satisfied guarantee} the asymptotic stability of the equilibria of the interconnected system (\ref{eq:network4}) and (\ref{eq:model}). We note here that these conditions are decentralized. \il{Before} presenting the aforementioned conditions, we first provide the following definition for input-strict passivity \cite{Kha01}.

\begin{definition} \label{definition:inputstrict}\il{
Consider \il{a dynamical} system represented by the state space model
\begin{equation} \label{eq:dynsys}
\begin{split}
\dot{x_i} =& f_i(x_i,u_i) \\
y_i =& g_i(x_i,u_i)
\end{split}
\end{equation}
where $f_i : \mathbb{R}^{n_i} \times \mathbb{R}^{p_i} \to \mathbb{R}^{n_i}$ and $g_i : \mathbb{R}^{n_i} \times \mathbb{R}^{p_i} \to \ill{\mathbb{R}^{p_i}}$ are locally \il{Lipschitz.}
% for any $n$ and $p \in \mathbb{N}$.
%Furthermore, the system (\ref{eq:dynsys}) has the same number of inputs and %outputs.
\il{Such system} is said to be locally input strictly passive \il{about the equilibrium $(\hat{u}_i, \hat{x}_i)$,} if there exist open neighborhoods $U_i$ and $X_i$ about \il{$\hat{u}_i,\ \hat{x}_i$, respectively,} \il{a continuously} differentiable function $\mathcal{V}_i(x_i)$ (called the storage function), \il{and a function $\phi(.)$} such that
\begin{equation} \label{eq:inputstrictly}
(u_i-\hat{u}_i)^\textrm{T} (y_i-\hat{y}_i) \geq \dot{\mathcal{V}_i} + (u_i-\hat{u}_i)^\textrm{T} \phi_i (u_i-\hat{u}_i)
\end{equation}
\il{for all $u_i \in U_i$ and all $x_i \in X_i$, where $(u_i-\hat{u}_i)^\textrm{T} \phi_i (u_i-\hat{u}_i) > 0$ for $u_i\neq\hat u_i$.}}
\end{definition}
\il{
\begin{remark}
For linear systems or systems linearized about an equilibrium the passivity property can be easily verified by means of computationally efficient methods using the KYP lemma \vs{\cite{kottenstette2014relationships}}, or via the positive realness of the transfer function. \vsp{We note here that \icl{the KYP lemma also allows to explicitly construct the storage function of the system, which for linear systems is a quadratic function of the form $\mathcal{V}(x_i)=x_i^{\textrm{T}} P_i x_i$, where matrix $P_i\in\mathbb{R}^{n_i\times n_i}$ is obtained by solving a convex optimization problem (a semidefinite program).  For nonlinear systems it can be verified by exploiting structural properties such as feedback interconnections of passive systems, or via an explicit construction of the storage function (see e.g.  \cite{kasis2016primary}, \cite{mccourt2013demonstrating}).}} \icl{An alternative way to check the passivity property for linear systems is via the positive realness of their transfer function $G(s)$. In particular, input strict passivity is implied if $G(j\omega) + G^\textrm{T} (-j\omega)$ is positive definite (or equivalently has positive eigenvalues) for all $\omega\in\mathbb{R}$.}
%This can be useful for control design, as it will be discussed in section %\ref{sec:simulations}.}
%\icl{????say something about positive realness and the eigenvalue condition in the frequency %domain????}
%\vsp{The aforementioned storage functions are usually quadratic and can be deduced using %Linear Matrix Inequalities (LMIs). Such an example is presented in our previous work %\cite{kasis2016primary}, where we constructed the storage functions of the adopted %second-order generator dynamics and the controllable loads.

\end{remark}}

\il{In the assumptions below $\hat x$ is an equilibrium point of the interconnected system, and $(\hat u_i,\hat x_i)$ are the corresponding constant inputs and states of the bus dynamics (\ref{eq:model}) at this point.}

\begin{assumption} \label{assumption:passivedyn2}
%There exist open neighborhoods $U_i^V$ of $(\hat{V}_{a,i},\hat{V}_{b,i})$ and %$X_i$ of $\hat{x}_i$, and a continuously differentiable, positive semidefinite %function $\mathcal{V}_i^B (x_i)$ such that
\il{%Consider an equilibrium $\hat{x}$ of the interconnected system.
For each $i\in\mathcal{N}$, each of the bus dynamical systems (\ref{eq:model})
%with input $u_i=[-I_{a,i}, -I_{b,i}]$ and output $y_i=[V_{a,i}, V_{b,i}]$
satisfies} a local input-strict passivity property \il{about $(\hat u_i,\hat x_i)$}, in the sense described in Definition \ref{definition:inputstrict}. %\il{about an equilibrium of the interconnected  system}.
\end{assumption}

Similarly to the approach presented \ill{in \cite{kasis2016primary}, \cite{devane2016primary},} we \illl{assume} that the aforementioned passivity property holds without specifying the precise form of the bus dynamics. \ill{This} will allow us to include \ill{in the} %into the elaborated
stability analysis a broad class of bus dynamics and a variety of frequency and/or voltage control mechanisms.
%\vsp{Moreover, despite the fact that the passivity condition is local it can be employed in %order to verify whether the bus dynamics are passive or not under various operating %conditions, such as contigency events. However, as we are about to discuss in Section %\ref{sec:simulations}, this paper focus on the areas of small-disturbance voltage and %frequency stability, as these are defined in \cite{kundur2004definition}.}

%\begin{remark}
%For linear or linearized systems of the form (\ref{eq:model}), passivity can be verified numerically using the KYP (Kalman-Yakuvbovich-Popov) lemma  via a semidefinite program, i.e. a convex optimization problem. Employing KYP lemma provides us the opportunity to derive system's storage function as well. For further reading regarding the KYP lemma and the passivity numerical verification consult \cite{Kha01, boyd2004convex}.
%\end{remark}

Finally, to guarantee convergence, we will require two additional conditions on the behavior of the interconnected system (\ref{eq:network4}) and (\ref{eq:model}). These conditions \il{will be used in the %are important for our
proof of the convergence} result in Theorem \ref{theorem:converge}.
\il{
\begin{assumption} \label{assumption:inputoutput}
Consider the  dynamics (\ref{eq:model}) at bus $i$. When $u_i(t)=\hat u_i \  \forall t$, then $\hat x_i$ is \illl{asymptotically} stable, i.e. there exists neighbourhood $\tilde X_i$ about $\hat{x}_i$ s.t. for all $x_i(0)\in\tilde X_i$, we have $x_i(t)\rightarrow\hat x_i$ as~$t\rightarrow\infty$.
%There \il{exists} an open neighborhood $U_i^I$ of $(-\hat{I}_{a,i},\ -\hat{I}_{b,i})$ such that \il{for} any constant input $(-I_{a,i},\ -I_{b,i}) \ \in \ U_i^I$,  \il{the} bus dynamics (\ref{eq:model}) \il{asymptotically} converge to $\hat{x}_i \in X_i$ and $(\hat{V}_{a,i}, \hat{V}_{b,i}) \in U_i^V$.
\end{assumption}
\begin{remark}
Note that this condition is trivially satisfied in many cases as generation dynamics are usually open loop stable. The condition could also be relaxed to allow for integrators \ill{at some buses} (used in e.g, secondary control), but this is not done here for simplicity in the presentation.
\end{remark}}

\begin{assumption} \label{assumption:localminima}
The storage \il{functions $\mathcal{V}_i$ in Assumption \ref{assumption:passivedyn2}} have a strict local \il{minimum} at the point
\il{$\hat{x_i}$}.
%$(\hat{I}_{a,i} , \hat{I}_{b,i}, \hat{x}_i, \hat{V}_{a,i}, \hat{V}_{b,i})$.
\end{assumption}
\il{
\begin{remark}
This is a technical condition often satisfied. This is satisfied, for example, for any linear system if the latter is observable and controllable.
\end{remark}
}

%%%%
\subsection{Stability \il{Result}}

The passivity properties presented for both the network and the bus model are now exploited in order to show that the \il{equilibria %(\ref{eq:equilibrium})
of the system} (\ref{eq:network4}) and (\ref{eq:model}), for which the assumptions stated are satisfied, are asymptotically attracting. This \il{is stated in} the following theorem.

\begin{theorem} \label{theorem:converge}
\illl{Suppose %Assumptions \ref{assumption:passivedyn2} - \ref{assumption:localminima} hold and
there exists an equilibrium of the interconnected system (\ref{eq:model}), (\ref{eq:network4})
%and a neighborhood $S$ of the initial conditions and  where the
 for which the bus dynamics (\ref{eq:model}) satisfy Assumptions \ref{assumption:passivedyn2} - \ref{assumption:localminima} %\ref{assumption:passivedyn2} and \ref{assumption:inputoutput}
for all $i \in \mathcal{N}$.  Then this equilibrium is asymptotically stable, i.e. there exists an open neighbourhood $S$ about this point such that for all initial conditions $x(0) \in S$, the solutions of the system %interconnected system (\ref{eq:model}), (\ref{eq:network4})
converge to this point.}
%there exists an open set $S$ of initial conditions such that when  initial conditions $x(0) \in %S$, the solutions of the system (\ref{eq:model}), (\ref{eq:network4}) converge to this 5 %equilibrium point.}
%the equilibria defined in Definition \ref{definition:equilibria}.
\end{theorem}

\begin{remark}
%As it becomes clear from the above results, one of the main contributions of the %current paper is
\il{It should be noted that the} \illl{stability conditions in} Theorem \ref{theorem:converge} are decentralized \il{as they are conditions on the local bus dynamics.}
%and do not require any information neither for the power network topology nor %from the neighboring buses. \il{In particular,} we guarantee the asymptotic %stability of the equilibria by imposing \il{on the} bus dynamics certain %decentralized passivity \il{conditions.}
%which could be crucial for the future operation of the existing power systems %since are beneficial for the design of highly distributed control schemes.
\il{A distinctive feature of those is that the bus dynamics are formulated at the system reference frame, thus allowing to consider networks with losses as was discussed in Section \ref{sec:NetworkModel}.
%The conservatism of these conditions is to be tested in
In Section \ref{sec:simulations} it will be shown that these stability conditions are not conservative by applying those to real power networks with realistic data.}
\end{remark}

%%%%%%%%%%%%%%%%%%%%%%%%%%%%%%%%%%%%%%%%%%%%%%%%%%%%%%%%%%%%%%%%%%%%%%%%%%%%%%%%

\section{DISCUSSION} \label{sec:discuss}

%%%%
\subsection{Network}

As we mentioned before, the main difference between the proposed approach and the recent literature is that the analysis is carried out in the system reference-frame instead of each local machine reference-frame. \icl{It should be noted that even though this change of reference frame does not have an effect in a centralized stability analysis, it provides important benefits when stability is deduced by means of decentralized conditions. In particular, as it has been shown in the paper, a local passivity property at the bus dynamics in this reference frame is sufficient to deduce stability in a general network, without having to resort to a lossless assumption on the transmission lines.
%not make any difference when the stability analysis is carried out in a centralized manner, %due to the fact that all the necessary information regarding the power grid is available, %it was proved to be helpful when seeking to deduce stability by means of local conditions. %Particularly, it allows the consideration of more relaxed conditions for the power network %since by adopting an input/output formulation, we showed that the equations describing the %interconnections in any power grid with arbitrary topology constitute a passive system %without neglecting the network's lossy nature.}
It should also be noted that both the active and the reactive power flows across the network are taken into account in the models} and the bus voltage magnitudes are not considered to remain constant when sudden generation or load disturbances appear across the \il{power grid.}
%Additionally, the proposed framework does not require the adoption of %Kron-reduced models, \il{while preserving} the accuracy of the analysis.

In addition, the transformation into the system reference-frame %gave the opportunity to deal with
\illl{leads to}
simpler network equations. Specifically, by comparing \illl{equations} (\ref{eq:networkcomps}) with equations (\ref{eq:networkcompsdq}), it is easy to discern the complexity added when the analysis is carried out in the local machine reference-frame due to the existence of the sinusoids. This was the main reason why several simplifications were considered \il{in the stability analysis for}  power networks within the recent literature. \vsp{Such an important simplification} was the adoption of lossless power networks where equations (\ref{eq:networkcompsdq}) are reduced to
\begin{equation}\label{eq:networkcompsdq2}
\begin{split}
I_{q,i} =& \sum_{j=1}^{|\mathcal{N}|} \Big[ B_{ij} V_{q,j} \sin(\eta_{ij}) - B_{ij} V_{d,j}\cos(\eta_{ij}) \Big] \\
I_{d,i} =& \sum_{j=1}^{|\mathcal{N}|} \Big[ B_{ij} V_{q,j} \cos(\eta_{ij})+ B_{ij} V_{d,j} \sin(\eta_{ij}) \Big].
\end{split}
\end{equation}

Commonly, the previous equations are further simplified by assuming that $V_{d,i}=0 \ \forall \ i\in \mathcal{N}$, which leads to the following simpler and less accurate form:
\begin{equation}\label{eq:networkcompsdq3}
\small{ I_{q,i} = \sum_{j=1}^{|\mathcal{N}|} \Big[ B_{ij} V_{q,j} \sin(\eta_{ij}) \Big] \mbox{ and } I_{d,i} = \sum_{j=1}^{|\mathcal{N}|} \Big[ B_{ij} V_{q,j} \cos(\eta_{ij})\Big]}.
\end{equation}

In order to avoid counting in the generator's transient reactances and thus to facilitate the analysis, several works such as \cite{kasis2015stability, trip2016internal, zhao2014design, machowski2000decentralized} also considered that \vs{the $q$-axis bus voltage is equal to the $q$-axis transient emf, i.e. $V_{q,i}=E'_{q,i}$} \illl{(see also section \ref{sec:BusDynamics})}.
\todoiny{need to say what $E_{q,i}$ is, referring e.g. to the Machine models in the next section}

%%%%
\subsection{Bus Dynamics}\label{sec:BusDynamics}

The adoption \illl{of a broad} class of systems to represent bus dynamics
\illl{provides another advantage} of the proposed approach, since it
\illl{allows} \illl{to include dynamics associated with}
%the adoption of bus
%dynamics of
a variety of power system components such as generators, loads, inverter-based RES, \illl{and FACTS devices.} It also \illl{gives} the opportunity to use \illl{more accurate} higher order dynamical models and incorporate %either several voltage and/or
\il{voltage and frequency} control mechanisms at the same time. \illl{Dynamic models of the synchronous generator are commonly used %found \illl{used}
in power system stability studies,} %dynamics.
\illl{ as %This is due to the fact that
their} stable operation guarantees the security and the reliability of a power system. \illl{The} adoption of more accurate generator models \illl{is hence} crucial, since \illl{simpler} ones \illl{can} fail to accurately predict the behavior of the system~\cite{venezian2016warning}.

We \ill{describe} below the fourth-order generator model which is widely considered to be sufficiently accurate to analyze electromechanical dynamics, and present how this can be incorporated into our approach. The aforementioned model which will be also used to verify our results in Section \ref{sec:simulations}, is described by the following set of differential equations:
\begin{equation} \label{eq:fourthorder}
\begin{split}
M_i \Delta \dot{\omega}_i =& P_i^m - P_i^e - D_i \Delta \omega_i \\
\dot{\delta}_i =& \Delta \omega_i \\
T'_{do,i} \dot{E'}_{q,i} =& E_{f,i} - E'_{q,i} + I_{d,i} (X_{d,i} - X'_{d,i})  \\
T'_{qo,i} \dot{E'}_{d,i} =& - E'_{d,i} - I_{q,i} (X_{q,i} - X'_{q,i})
\end{split}
\end{equation}
where the electrical power $P_i^e =  E'_{q,i} I_{q,i} + E'_{d,i} I_{d,i} + (X'_{d,i} - X'_{q,i}) \ I_{d,i} \ I_{q,i}$. The synchronous generator described in (\ref{eq:fourthorder}), is modeled in its local $dq$-reference frame and it is represented by the transient emfs $E'_{d,i}$ and $E'_{q,i}$ behind the transient reactances $X'_{d,i}$ and $X'_{q,i}$ as defined by the following equation
\vsp{\begin{equation} \label{eq:dqvoltage}
\begin{bmatrix} V_{q,i} \\ V_{d,i} \end{bmatrix} = \begin{bmatrix} E'_{q,i} \\ E'_{d,i} \end{bmatrix} - \begin{bmatrix} R_i & -X'_{d,i} \\  X'_{q,i} & R_i \end{bmatrix} \begin{bmatrix} I_{q,i} \\ I_{d,i} \end{bmatrix}.
\end{equation}}

\illl{In particular,} the synchronous generator model forms a 2-input/2-output system with inputs the currents $-I_{d,i}$ and $-I_{q,i}$, and outputs the voltages $V_{d,i}$ and $V_{q,i}$. In order to allow the coupling with the network model (\ref{eq:network4}), the generator dynamics have to be transformed into the system reference frame. We therefore incorporate the mappings (\ref{eq:mapping}) - (\ref{eq:invmapping2}) into \illl{the} synchronous generator dynamics, and the current $dq$-components $I_{d,i}$ and $I_{q,i}$ are replaced in equations (\ref{eq:fourthorder}) by the net injected current components $I_{a,i}$ and $I_{b,i}$ using (\ref{eq:mapping2}). The outputs $V_{a,i}$ and $V_{b,i}$ are derived by (\ref{eq:dqvoltage}) as follows
\vsp{\begin{equation} \label{eq:abvoltage}
\begin{bmatrix} V_{a,i} \\ V_{b,i} \end{bmatrix} = T_i^{-1} \begin{bmatrix} E'_{q,i} \\ E'_{d,i} \end{bmatrix} - T_i^{-1} \begin{bmatrix} R_i & -X'_{d,i} \\  X'_{q,i} & R_i \end{bmatrix} T_i \begin{bmatrix} I_{a,i} \\ I_{b,i} \end{bmatrix}
\end{equation}}
\begin{table}[ht]
\centering
\caption{Description of the Main Variables and Parameters appearing in the Synchronous Generator Models}
\label{description}
\footnotesize{
\vspace{-.2cm}
\begin{tabular}{ll|ll}
\multicolumn{2}{c|}{\textbf{Variables}} & \multicolumn{2}{c}{\textbf{Parameters}} \\ \hline
$\Delta \omega$  & Frequency deviation              &  $M$     & Moment of inertia                   \\
$\delta$         & Stator's phase angle             &  $D$     & Damping coefficient                 \\
$E'_d $          & d-axis transient EMF             &  $X_d $  & d-axis synchronous reactance        \\
$E'_q$           & q-axis transient EMF             &  $X'_d$  & d-axis transient reactance          \\
$I_d$            & d-axis current    							  &  $T'_d$  & d-axis open-circuit time constant   \\
$I_q$            & q-axis current							      &  $X_q$   & q-axis synchronous reactance        \\
$V_d$            & d-axis bus voltage               &  $X'_q$  & q-axis transient reactance          \\
$V_q$            & q-axis bus voltage               &  $T'_q$  & q-axis open-circuit time constant   \\
$P^m$            & Mechanical power                 &  $R_i$   & stator windings resistance          \\
$P^e$            & Electrical power                 &          &\\
$E_f$            & exciter output emf               &          &\\
$E_g$            & controlled voltage               &          &
\end{tabular}}
\end{table}

We can now identify that the aforementioned generator model matches to the class of bus models (\ref{eq:model}). The dynamic model of the generator corresponds to the vector function $f_i$ while equations (\ref{eq:abvoltage}) match to the vector functions $g_i$. The presented formulation still holds for the \illl{second} and the third-order synchronous generator models, where %We only need to consider that
the transient emfs $E'_{d,i}$ and $E'_{q,i}$, respectively, are assumed to remain constant. Higher order models, such as the fifth or the sixth order models can be also incorporated in this framework in an analogous way. In these models, the synchronous generators are represented by their subtransient emfs $E^{''}_{d,i}$ and $E^{''}_{q,i}$ behind the subtransient reactances $X^{''}_{d,i}$ and $X^{''}_{q,i}$. More detailed information about generator modeling can be found in \cite{kundur1994power, machowski2011power, sauer1997power}.

\begin{figure}[!htbp]
\centering
\includegraphics[width=\columnwidth]{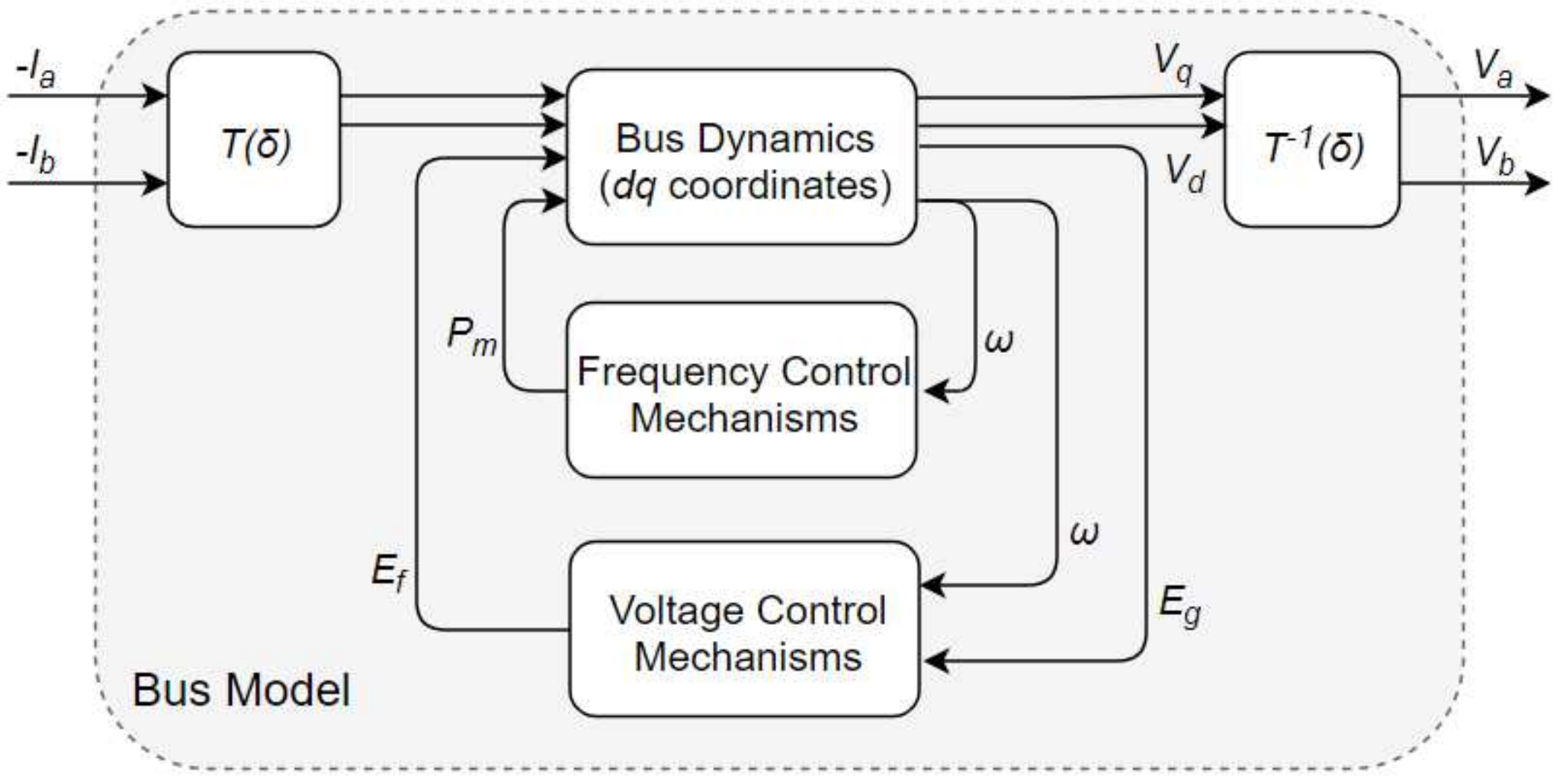}
\caption{A graphical representation of a generator model expressed in system reference frame.}
\label{fig:busmodel}
\end{figure}

As we discussed earlier, the generator dynamics can be expanded further so as to introduce the dynamics of \illl{frequency} and/or voltage control mechanisms, thus allowing the derivation of more accurate stability results. Such frequency and voltage control mechanisms include several types of turbine governors, exciters and \icl{power system stabilizers}. A graphical representation of the adaption of a synchronous generator model to the proposed framework along with the incorporation of frequency and voltage control is provided in Figure \ref{fig:busmodel}.
We highlight here the fact that the passivity conditions which ensure the asymptotic stability of the equilibria, refer to the bus dynamics and not specifically to the voltage and the frequency control systems that are applied. This allows us to include \ill{more advanced} regulation mechanisms which in \il{most cases} are not passive (e.g. turbine governors, excitation systems etc.). This \ill{is} an important advantage \illl{since} \ill{such dynamics are often omitted in approaches commonly presented in the related literature, due to the additional complexity they introduce.}
%in contrast to the common approaches presented in the related literature, the %incorporation of such dynamics makes analysis more \ill{involved and are often %omitted}. %complex and difficult.% For instance in \cite{stegink2016optimal}, %although T. Stegnik et al. managed to construct an optimal frequency controller %in power networks of high-dimensional models (sixth order) of synchronous %machines, they avoided to adopt any voltage control mechanism due to the %complexity that would be added to their analysis.

Finally, \il{it should be noted that loads, inverter-based RES and FACTS devices can also easily be incorporated within our approach}. However, \il{an explicit analysis of such models is beyond the scope} of the current paper and is therefore omitted. Several detailed dynamic models for the aforementioned power system components can be found within \cite{kundur1994power, hiskens1995load, schiffer2016survey, teodorescu2011grid, wang1997unified, wang1998unified, wang2000unified}. It is crucial to mention here that the framework becomes more flexible and easy to apply when we deal with bus dynamics that do not \il{involve a} rotating axis. In such cases, bus dynamics are already expressed in \illl{the} system reference frame and thus there is no need to include the maps $T_i$ and $T_i^{-1}$ into the adopted dynamics.

%%%%%%%%%%%%%%%%%%%%%%%%%%%%%%%%%%%%%%%%%%%%%%%%%%%%%%%%%%%%%%%%%%%%%%%%%%%%%%%%

\section{SIMULATIONS} \label{sec:simulations}

In this section we verify our framework and the derived stability results through \ill{applications} on the Two Area Kundur Test System \cite{kundur1994power} and the IEEE New York / New England 68-bus interconnection system \cite{rogers2012power}. These applications focus \ill{on generator} buses and are carried out using the Power System Toolbox (PST) \cite{chow2000power}. Within the simulations, the generators are modeled by the fourth-order dynamics (\ref{eq:fourthorder}) on which both frequency and voltage control mechanisms are applied. Specifically, \ill{frequency and voltage control} are carried out by turbine governors and exciters respectively, while PSSs are applied to the generator's excitation system. The adopted models of the turbine governors, the exciters and the PSSs can be found in PST manual \cite{chow2000power}.

%are described in Laplace domain by
%\begin{equation*}
%P_i^m = \frac{(1+sT_{3,i})(1+sT_{4,i})}{(1+sT_{s,i})(1+sT_{c,i})(1+sT_{5,i})} (P_{ref,i}^m - K_i^{TG} \Delta \omega_i )
%\end{equation*}
%and
%\begin{equation*}
%E_{f,i} = \frac{K_i^A}{(1+sT_{A,i})} \Big( V_i^{ref} + V_i^{PSS} - \frac{1}{(1+sT_{R,i})} |E_{g,i}|\Big)
%\end{equation*}
%respectively. For the PSSs we use the conventional PSS model which is also given in Laplace domain by
%\begin{equation*}
%V_i^{PSS} = K_i^{PSS} \frac{sT_{w,i} (1+sT_{1,i})(1+sT_{3,i})}{(1+sT_{w,i})(1+sT_{2,i})(1+sT_{4,i})} \Delta \omega_i.
%\end{equation*}
%In the above transfer functions, the variables defined with the letter T and the letter K denote the time constants and the gains of the respective control mechanism. Additionally, $P_{ref,i}^m$ is the mechanical power reference input of the turbine governor, $V_i^{ref}$ is the voltage reference input of the exciter and $V_i^{PSS}$ is the supplementary input injected by the PSS to the exciter. Finally, $|E_{g,i}|$ represents the magnitude of the generator's terminal voltage and is given by $|E_{g,i}| = \sqrt{E_{q,i}^2 + E_{q,i}^2}$.

\vsp{In order to facilitate the verification of the passivity property on the generator buses of the test systems, we linearize the dynamics of each generator bus individually about an \icl{equilibrium}.}
%We should mention here that the linearization is carried out only for each generator %dynamics and not for all the power system components as in Small Signal Stability Analysis %\cite{kundur1994power}.}
The equilibria are identified by solving a \illl{Power Flow} problem for each test system \illl{respectively\footnote{\illl{It should be noted that the phase difference $\delta_i$ between each local (d,q) and the system reference frame is obtained from each generator's q-axis transient emf $E'_{q,i}$, rather than the q-axis bus voltage $V_{q,i}$.}}}.
%For the derivation of the transfer matrices of the generator buses, we linearize the dynamics of each generator bus individually about an equilibrium. We should mention here that the linearization is carried out only for each generator dynamics and not for all the power system components as in Small Signal Stability Analysis \cite{kundur1994power}.} The equilibria are identified by solving a \illl{Power Flow} problem for each test system \illl{respectively\footnote{\illl{It should be noted that the phase difference $\delta_i$ between each local (d,q) and the system reference frame is obtained from each generator's q-axis transient emf $E'_{q,i}$, rather than the q-axis bus voltage $V_{q,i}$.}}}.
%, \vsp{at its initial conditions}
%\footnote{\vsp{The resulting equilibria when the test systems are subjected to a %significant load change, are very close to those derived at the initial conditions.}}
%\vsp{The phase difference $\delta_i$ between each local (d,q) and the system reference %frame is obtained from each generator's q-axis transient emf $E'_{q,i}$ and not the %q-axis bus voltage $V_{q,i}$.}
In order to verify the passivity of the bus models we use Linear Matrix Inequalities (LMIs) whose application on passivity verification is extensively described in \cite[Section~2]{mccourt2013demonstrating}. An alternative way to verify the passivity of bus dynamics with transfer matrix $G_i(s)$, is by checking the positive definiteness of the matrix \ill{$G_i(\textrm{j}\omega) + G_i^\textrm{T}(-\textrm{j}\omega)$} as indicated in \cite{kottenstette2010relationships}. In particular, the positive definiteness
%of each bus transfer function matrix
is ensured when the eigenvalues \icl{of the matrix are positive\footnote{\icl{Note that the eigenvalues of $G_i(\textrm{j}\omega) + G_i^\textrm{T}(-\textrm{j}\omega)$ are always real as the matrix is Hermitian.}}}.

We first deal with the Four Machine Two-Areas Kundur Test System which is widely used for stability studies. The passivity of the four generator buses is verified using LMIs, for the following four different cases: (i) no turbine governor / no exciter / no PSS, (ii) turbine governor / no exciter / no PSS, (iii) turbine governor / exciter / no PSS and (iv) turbine governor / exciter / PSS. All generator buses are not passive when neither of the available control mechanisms is employed. When turbine governors are added to generators, bus dynamics are slightly damped but they still remain non passive. The exciters further passivate the generator buses making buses 1 and 2 passive. Buses 11 and 12 still remain non passive. Finally, the application of PSSs to the generators completely passivates the dynamics.

\begin{figure}[t!]%[htbp!]
\centering
\includegraphics[width=\columnwidth]{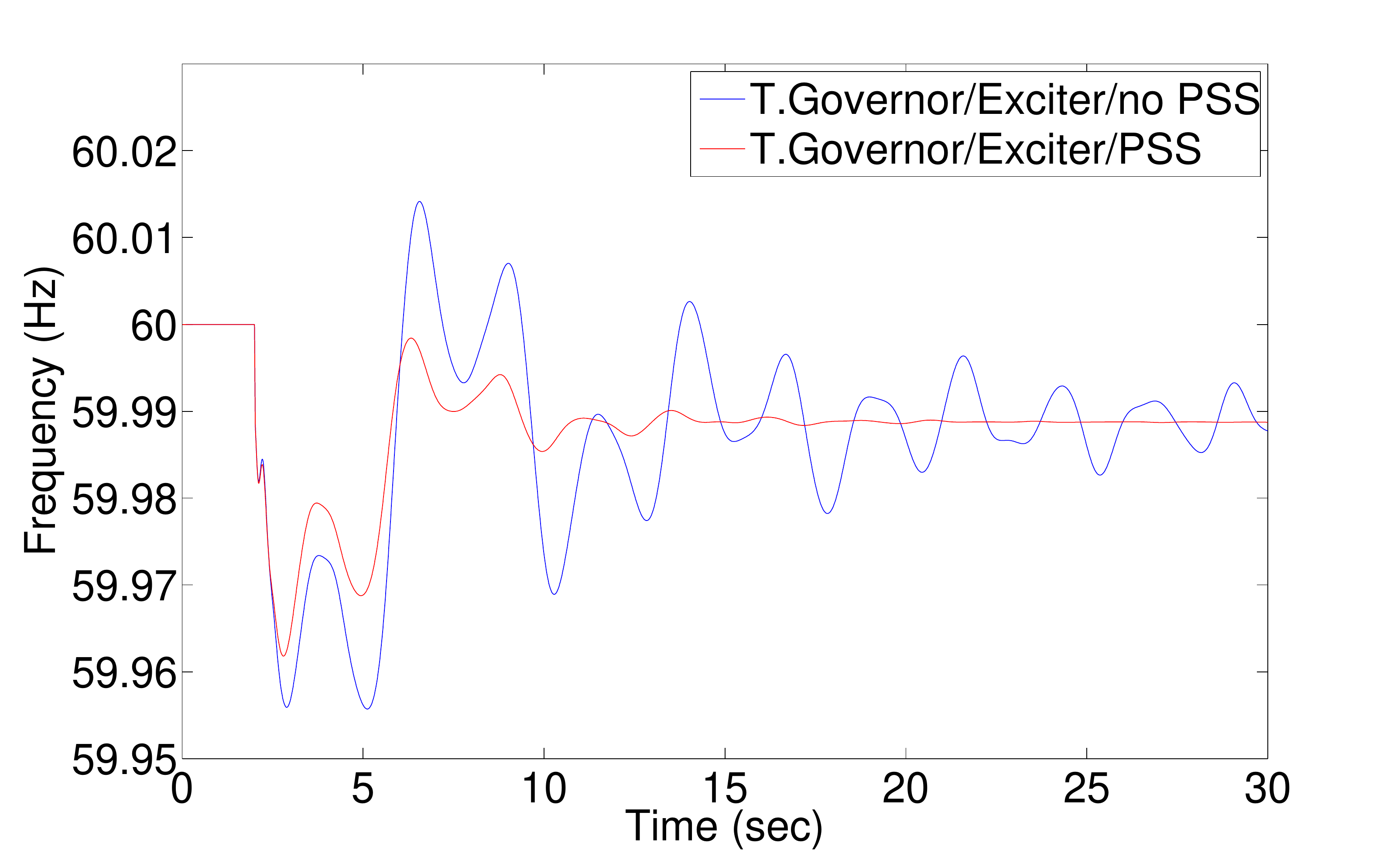}
\caption{Frequency deviation at bus 27 after a sudden change of 1pu at the load buses 1, 9 and 18.}
\label{fig:frequency}
\end{figure}

\begin{figure}[t!]%[htbp!]
\centering
\includegraphics[width=\columnwidth]{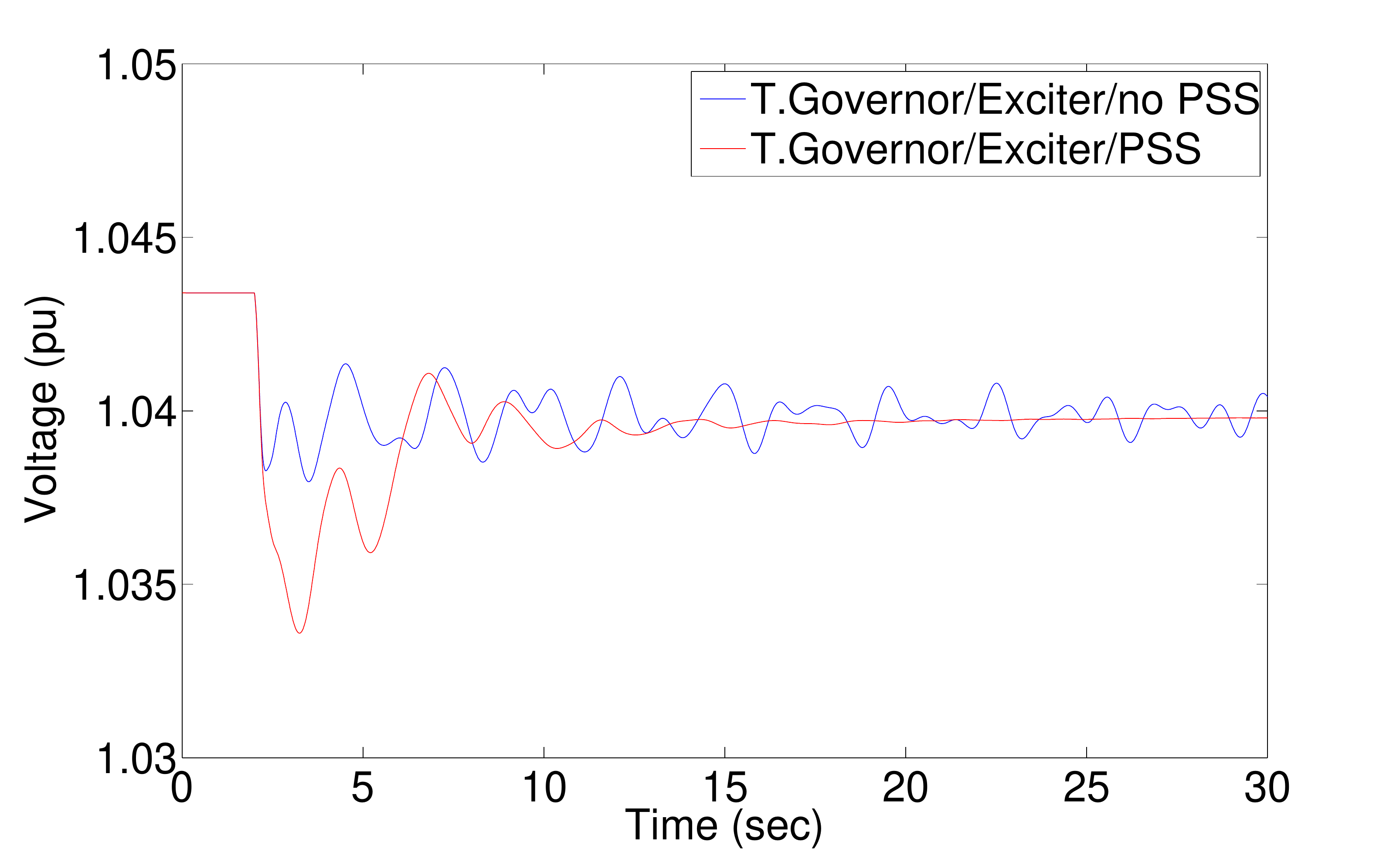}
\caption{Voltage deviation at bus 27 after a sudden change of 1pu at the load buses 1, 9 and 18.}
\label{fig:voltage}
\end{figure}

The proposed approach is also applied on the generator buses of the IEEE 68 bus test system.
%The tuning of the PSSs that are applied to the generators' exciters is derived in \cite{cai2003simultaneous}.
According to the derived results, the generator buses are also non passive for the cases (i) and (ii). In both cases, the power system collapses after a sudden change of load across the network. On the other hand, when the excitation system is applied to the generators, the generator buses are considerably damped and, although the system presents an oscillatory behavior, it remains stable when a generation-load mismatch occurs. To be more specific, the application of the excitation system on the grid-connected generators makes generator buses \vsp{53, 59, 61 and 64} passive while the rest remain non passive. Finally, the incorporation of the PSSs at the generator exciters passivates further the generator buses, and results \illl{in} a more stable and robust operation. All generator buses are now \illl{passive\footnote{\illl{It should be noted that the passivity property was verified for all choices of reference bus for the angles $\delta_i$. The choice of reference did not affect the passivity property since the relative values of the angles are close to 0.}}} except from buses \vsp{58, 62, 63 and 65}. However, we can achieve to \illl{passivate} these generator buses by slightly increasing the transient reactances $X'_d$ and $X'_q$ of the respective generators (approximately $15\%$). These results are illustrated in Figures \ref{fig:frequency} and \ref{fig:voltage} which present the frequency and the voltage deviation at bus 27 respectively, when a sudden change of 1pu is applied at the load buses 1, 9 and 18. \vsp{This change corresponds to a total change of 300 MW. We should note here that the IEEE 68 bus test system consists of a total load of \icl{18.33} GW.} Due to the fact that the system collapses when the excitation system is not applied to the generators, we omit the respective figures for the cases (i) and (ii). We should also note here that, although not all the generator buses are passive, the power system is stable\footnote{\icl{Passivity is a sufficient condition for stability
%a more strict property than stability
which implies that the power system can be stable even if not all buses are passive.
%Although passivity is a sufficient condition to guarantee the asymptotic stability of the %system, it cannot be also considered as necessary when used in a decentralized manner.
It should be noted that its essence is that it is a decentralized condition, and any decentralized stability condition is in general only sufficient as
in order to derive a necessary and sufficient stability condition, the explicit knowledge of the dynamics of the whole power grid is required.}}. From these results, it becomes clear that under a proper design of the control mechanisms which are applied to the generators, we can achieve to completely passivate the generator buses of the network, and thus to ensure the asymptotic stability of the power system.
%\vsp{Finally, it should be also mentioned that any arbitrary selection of reference, %does not affect the verified passivity property since the angles $d_i$ at generator %buses are very close to 0.}
%\ill{????? does the equilibrium point change significantly with the disturbance added? The passivity property was checked before the disturbance ?????}
%\ill{????Increase the fonts in the figures????}


\begin{figure}[t!]%[htbp!]
\centering
\includegraphics[width=\columnwidth]{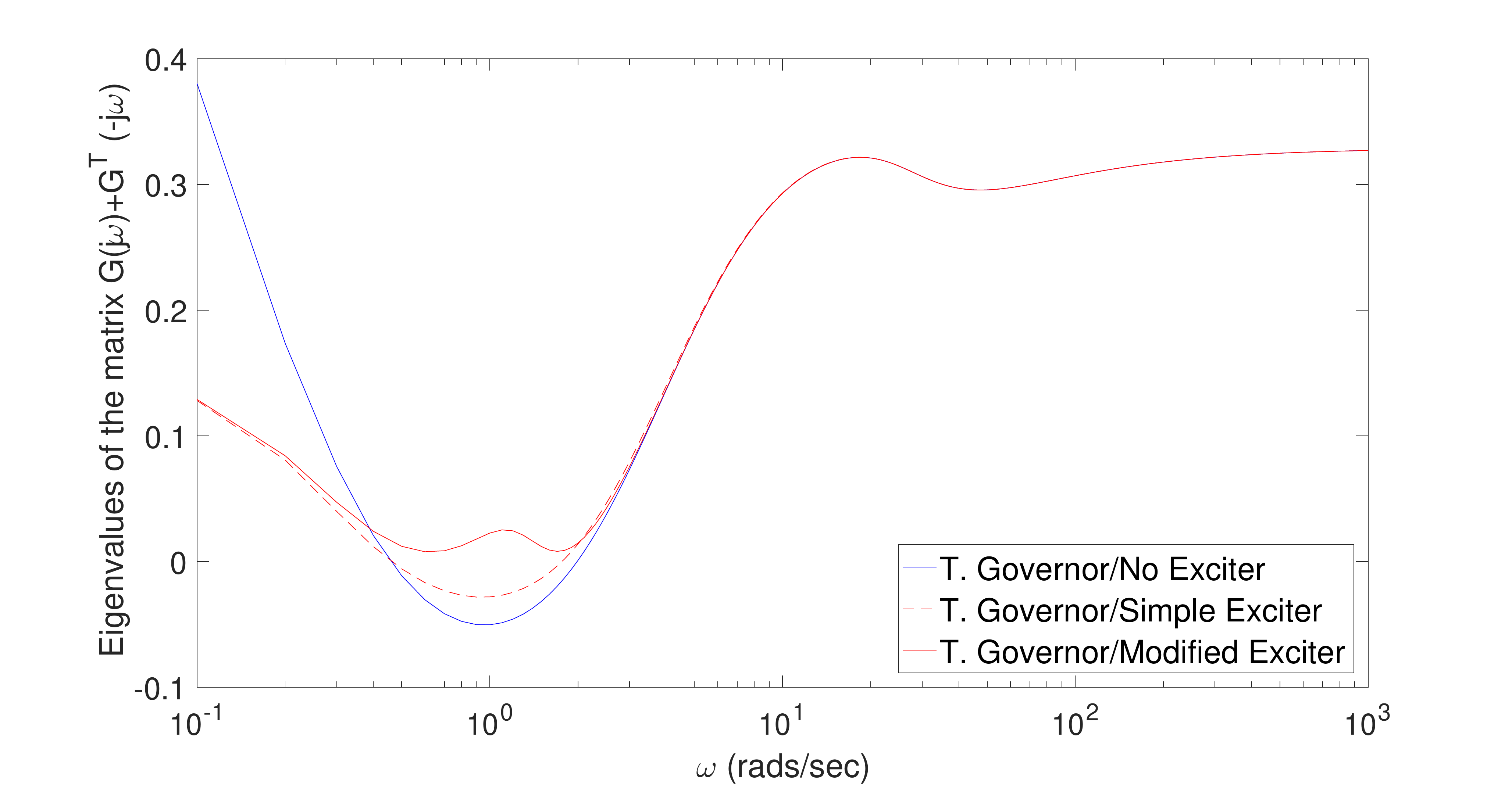}
\caption{Eigenvalues of $G(j\omega) + G^\textrm{T} (-j\omega)$, where $G(s)$ is the transfer function of the linearized dynamics \eqref{eq:model} at generator bus 53. The figure shows the eigenvalues in the problematic range where the passivity property is violated.}
\label{fig:loci}
\end{figure}

\vsp{\icl{{\color{black}In order to demonstrate how existing control mechanisms can be designed so as to satisfy the passivity property, we consider a generator bus with turbine governor where the bus dynamics are non passive without excitation control, and also a simple first order exciter leads to marginally non passive dynamics.
%We then discuss how incorporating a more advanced appropriately tuned lead/lag compensator %in the exciter
%we introduced a specific turbine governor setting that even when the exciters are %incorporated to the analysis the majority of bus dynamics remain at least marginally %non-passive.
We then discuss how modifing the exciter
%implementation, i.e.
%\begin{equation*}
%E_{f,i}=\frac{K_a}{1+sT_a} E_{g,i}
%\end{equation*}
by adding an additional phase lag compensator can passivate the dynamics.

In particular, %the first  order and
the modified exciter has transfer function given by
%\begin{equation*}
%E_{f,i}=\frac{K_a}{1+sT_a} E_{g,i}
%\end{equation*}
%and
\begin{equation}
E_{f,i}=\frac{K_a}{1+sT_a} \frac{1+sT_c}{1+sT_b} E_{g,i}
\label{eq:exciter}
\end{equation}
where the term $(1+sT_c)/(1+sT_b)$ is the phase lag compensator that has been added such that the passivity property is satisfied. We describe below in detail how the parameters of the transfer function in \eqref{eq:exciter} have been chosen and provide their values in Table \ref{values} below.
%The values chosen for the constants in the transfer functions above are given in the table %provided below. %We should note here that the constants $K_a$ and $T_a$ which are used within both the %simple and the modified exciter, are common.
\begin{table}[htbp!]
\centering
\caption{\color{black}Values of the constants of the modified exciters.}
\label{values}
\begin{tabular}{ccccc}
\hline
\textbf{Generator}    & $\bm{K_a}$    & $\bm{T_a}$    & $\bm{T_b}$    & $\bm{T_c}$       \\ \hline
53                    & 20             & 0.05           & 3.3            & 0.3           \\
54                    & 20             & 0.05           & 1.5            & 0.1           \\
55                    & 20             & 0.05           & 2.5            & 0.2           \\
56                    & 20             & 0.05           & 1.5            & 0.5           \\
57                    & 20             & 0.05           & 1.7            & 0.3           \\
58                    & 20             & 0.05           & 1.4            & 0.7           \\
59                    & 20             & 0.05           & 1.5            & 0.4           \\
60                    & 20             & 0.05           & 1.8            & 0.2           \\
61                    & 20             & 0.05           & 1.9            & 0.1           \\
62                    & 20             & 0.05           & 1.8            & 0.4           \\
63                    & 20             & 0.05           & 1.0            & 0.1           \\
$64^{**}$             & 20             & 0.05           & 0.0            & 0.0           \\
$65^{**}$             & 20             & 0.05           & 0.0            & 0.0           \\
$66^{**}$             & 20             & 0.05           & 0.0            & 0.0           \\
$67^{**}$             & 20             & 0.05           & 0.0            & 0.0           \\
$68^{**}$             & 20             & 0.05           & 0.0            & 0.0           \\ \hline
%\multicolumn{5}{l}{* : Generator remains passive under any selection of governor %parameters} \\
\multicolumn{5}{l}{** : No modifications were applied to the generator.}                      \\ \hline
\end{tabular}
\end{table}

The gain $K_a$ is the dc gain of the exciter and the parameter $T_a$ determines the cutoff frequency (given by $1/T_a$) of the first order system $1/(1+sT_a)$. Very large values of $K_a$ and $1/T_a$ will violate the passivity of the bus dynamics and lead to oscillatory responses. On the other hand, small values of $K_a$ and $1/T_a$ will lead to a feedback scheme that is slow in its response. Therefore  $K_a$ and $1/T_a$ are chosen large enough while also ensuring that the passivity property is not violated.

The selection of the values of the time constants $T_b$ and $T_c$  of the phase lag compensator is carried out so that a reduction in gain is achieved in the problematic frequency range where the passivity property is violated, i.e. the range of frequencies where the eigenvalues of $G(j\omega) + G^T (-j\omega)$ are negative, where $G(s)$ is the transfer matrix of the corresponding bus dynamics. The phase lag compensator is therefore designed such that $1/T_b$ is smaller than the problematic frequencies. Parameter $T_c$ is then chosen such that $T_c<T_b$ with $T_c$ sufficiently small so as to achieve a sufficient reduction in gain, and also $1/T_c<1/T_a$ so as to avoid additional phase lag at higher frequencies.

This process is demonstrated in Figure \ref{fig:loci} where the eigenvalues of the matrix $G(j\omega) + G^\textrm{T} (-j\omega)$ are illustrated at different frequencies, where $G(s)$ is the transfer function of the linearized dynamics \eqref{eq:model} at generator bus 53. The figure shows the eigenvalues\footnote{\icl{Note that $G(s)$ is a $2\times2$ matrix and for convenience in the presentation only one of the two eigenvalues is shown where the passivity condition is violated.}} in the regime where the passivity property is violated. In particular, as seen from the figure the passivity property is violated in the frequency range {\color{black}$\omega \in [0.3, 3]$} rads/sec  (since the eigenvalues are negative) when no exciter or the simple first order exciter are used. The figure also shows that the addition of an appropriately tuned phase lag compensator passivates the dynamics.} It should be noted that the essence of this design process is that it is decentralized, based on only the local bus dynamics, without requiring at each bus to be aware of the dynamics of the entire network as in a classical small signal analysis.}}

%The passivation of the dynamics of generator bus 53 is presented in Figure \ref{fig:loci} through the graphical illustration of the loci of the real part of the problematic eigenvalue for $\omega \in [0.1, 1000]$ rads/sec.}

\begin{figure}[t!]%[htbp!]
\centering
\includegraphics[width=\columnwidth]{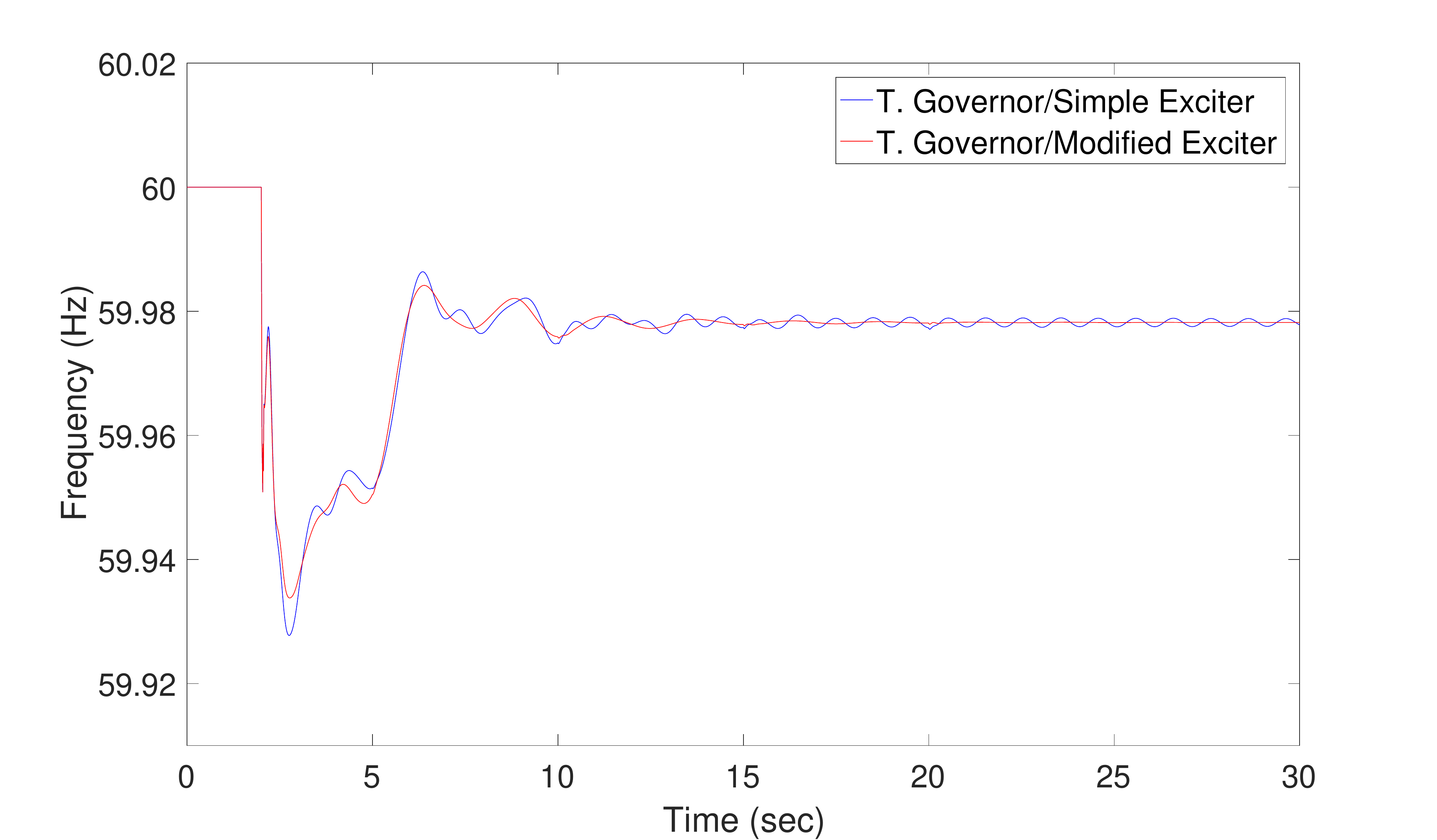}
\caption{Frequency deviation at bus 27 after a sudden change of 3pu at the load buses 1, 9, 18, 20, 37 and 42.}
\label{fig:frequency2}
\end{figure}

\begin{figure}[t!]%[htbp!]
\centering
\includegraphics[width=\columnwidth]{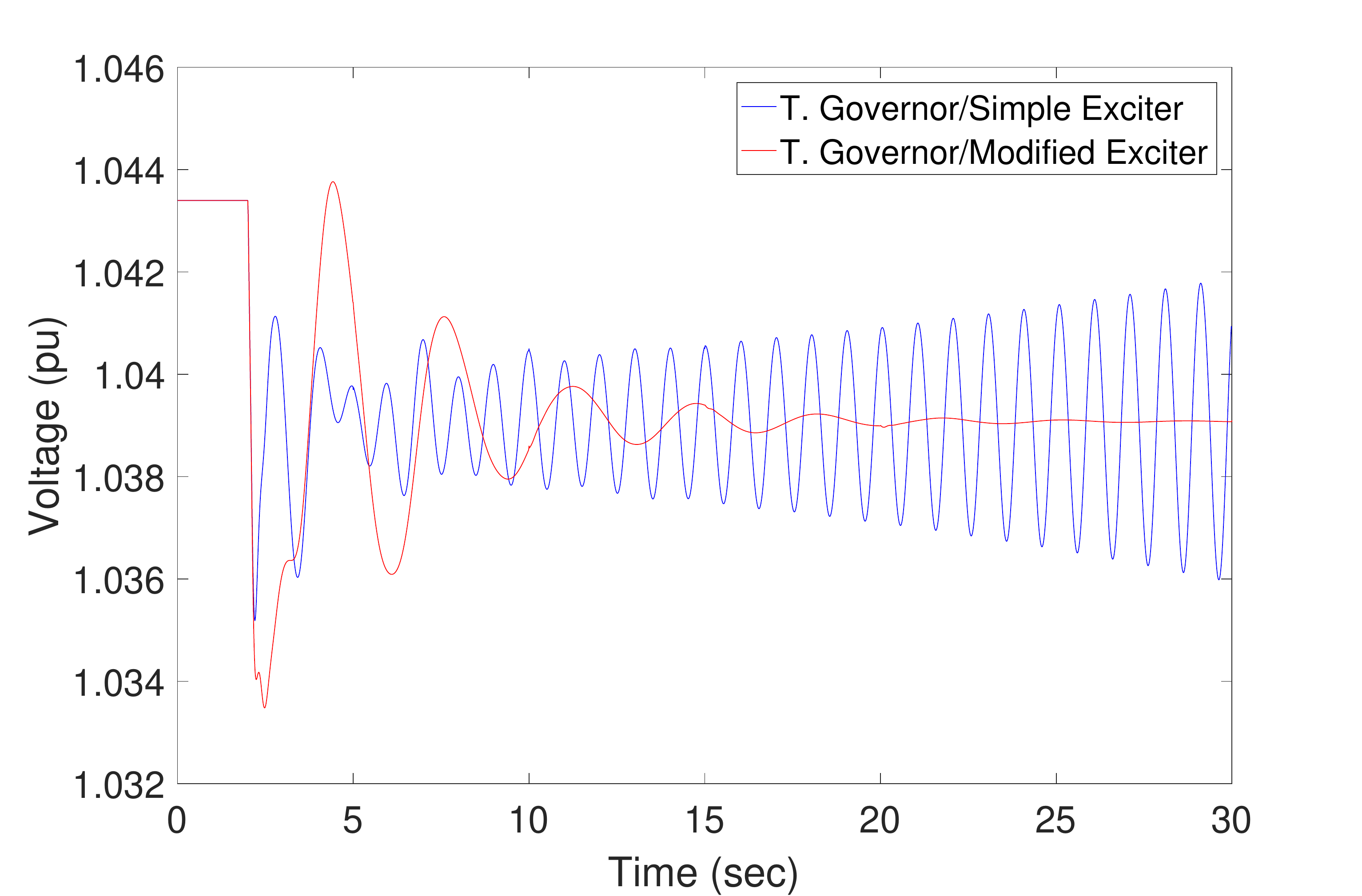}
\caption{Voltage deviation at bus 27 after a sudden change of 3pu at the load buses 1, 9, 18, 20, 37 and 42.}
\label{fig:voltage2}
\end{figure}

\vsp{The \icl{performance of the system when the passivity based design described above is applied at all buses within the network is illustrated} in Figures \ref{fig:frequency2} and \ref{fig:voltage2}. \icl{The figures present} the frequency and the voltage \icl{deviation, respectively,} at bus 27, when a sudden change of 3pu is applied at the load buses 1, 9, 18, 20, 37 and 42. As it can be seen from both the aforementioned figures, the introduction of an appropriately tuned lag \icl{compensator} to the excitation system of the generator results in a significantly less oscillatory behavior of the system. We should note here that for the \icl{dynamic} simulation, we considered a total load change of $1800$ MW which correspond to $10\%$ of the grid-connected \icl{load}. }
%This implies that the passivity property of the generator buses is satisfied for relatively %large deviations from the initial equilibria despite the fact that the condition presented in Assumption \ref{assumption:passivedyn2} is local.}

\begin{figure}[t!]%[htbp!]
\centering
\includegraphics[width=\columnwidth]{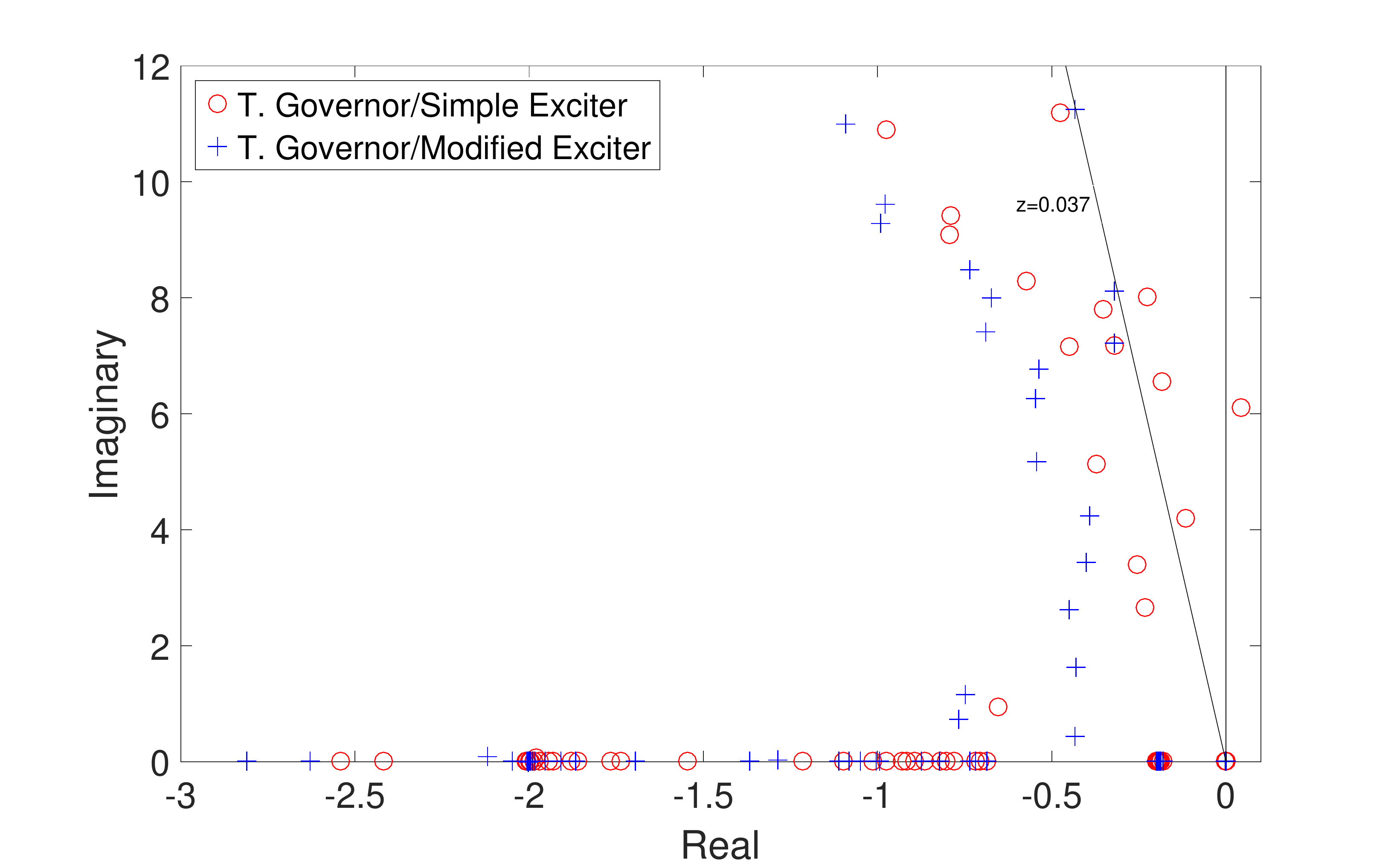}
\caption{Eigenvalues of the linearized network dynamics of the IEEE 68 bus test system.}
\label{fig:small}
\end{figure}

\vsp{The stability enhancement achieved through the modification of the initial simple exciter is also illustrated through the eigenanalysis of the test system. As it can be seen from \icl{Figure~\ref{fig:small}, the application of a lag compensator to the exciter that passivates the bus dynamics} significantly damps the calculated modes. More specifically, when the simple exciter \icl{that violates the passivity property} is employed to the generators the system is small-signal unstable since there exists an eigenvalue with positive real part. On the other hand, the application of the modified exciter on the generators stabilizes the system moving all the eigenvalues to the left half plane. Moreover, the proposed \icl{compensator} to the excitation system of the generators yields \icl{a good damping ratio\footnote{\vsp{The damping ratio of an eigenvalue $\lambda = a + j \beta$ is defined as $z = - \frac{a}{\sqrt[]{a^2 +\beta^2}}$.}} for the modes} of the system.}

%%%%%%%%%%%%%%%%%%%%%%%%%%%%%%%%%%%%%%%%%%%%%%%%%%%%%%%%%%%%%%%%%%%%%%%%%%%%%%%%

%\section{CONCLUSIONS AND FUTURE WORK} \label{sec:conclusions}
\section{CONCLUSIONS} \label{sec:conclusions}

%%%%
%\subsection{Summary}
Within the paper, we have presented a novel passivity-based approach for decentralized stability analysis and control of power grids through the transformation of both the network and the bus dynamics into \illl{the} system reference frame. In particular, by adopting
%a $(2 \times |\mathcal{N}|)$-input/$(2 \times |\mathcal{N}|)$-output
\illl{an input/output}
formulation, the power networks were shown to constitute a passive system even if \ill{the} network's lossy nature is taken into account. We have then introduced a broad class of bus dynamics that fits to the aforementioned network formulation and provided \ill{passivity} conditions which when satisfied, \il{guarantee} the asymptotic stability of the \illl{entire} power system in a \il{decentralized} manner. Furthermore, we discussed \ill{the} opportunities \ill{and} advantages offered by this approach while explaining how the proposed framework can allow the inclusion of various power system components, such as synchronous generators \illl{and dynamic loads.} Finally, our approach \illl{has been} verified through \illl{various applications}
% involving}
%fourth-order synchronous generator models
\illl{in} the Two Area Kundur and the IEEE 68 bus test systems.

\section*{APPENDIX}

\begin{proofl}
By substituting the network equations (\ref{eq:network4}) in inequality (\ref{eq:condition1}) we get
\begin{equation} \label{eq:condition2}
\small{
\begin{split}
u^\textrm{T} y =& [V_a^\textrm{T} \ V_b^\textrm{T}]  H_{2n} \begin{bmatrix} V_a \\ V_b \end{bmatrix} =
[V_a^\textrm{T} \ V_b^\textrm{T}]  \begin{bmatrix} G_n & -B_n \\ B_n & G_n \end{bmatrix} \begin{bmatrix} V_a \\ V_b \end{bmatrix} \\
=& V_a^\textrm{T} G_n V_a + V_b^\textrm{T} G_n V_b \geq 0
\end{split}}
\end{equation}
for all $V_a$, $V_b$ $\in \ill{\mathbb{R}^{|\mathcal{N}|}}$. The inequality (\ref{eq:condition2}) reveals that the passivity of the network is ensured when the composite matrix $H_{2n}$ or equivalently its diagonal elements $G_n$, are positive semidefinite matrices.

$G_n \in \mathbb{R}^{|\mathcal{N}| \times |\mathcal{N}|}$ is a square, sparse symmetric matrix with non negative diagonal and negative off-diagonal elements, i.e., $G_{ii} \geq 0$ and $G_{ij} \leq 0$ $\forall \ i, \ j=1,2,\dots,|\mathcal{N}|$. It is \il{also diagonally dominant as the following equation holds:
\begin{equation}\label{eq:diagonaldom}
\small{
 G_{ii} = -\sum_{j\neq i}^{|\mathcal{N}|} G_{ij}  \Rightarrow |G_{ii}| = \sum_{j\neq i}^{|\mathcal{N}|} |G_{ij}|} \quad (i,j) \in \mathcal{E}
\end{equation}}
In order to prove the positive semidefiniteness of the matrix $G_n$, we now define the Geshgorin discs $D_i(G_{ii},R_i), \ i=1,2,\dots,|\mathcal{N}|$. $D_i$ is a closed disc centered at $(G_{ii},0)$, with radius $R_i=\sum_{i \neq j} |G_{ij}|$.
%\begin{equation*}
%R_i=\sum_{i \neq j} |G_{ij}|.
%\end{equation*}
As stated above the matrix $G_n$ has positive diagonal elements and is also \il{diagonally} dominant. Subsequently, its Geshgorin discs lie in the right half plane, have center on the real axis and are tangent to the imaginary axis since $G_{ii} - R_i = 0, \ \forall \ i=1,2,\dots,|\mathcal{N}|$. According to the Geshgorin circle theorem \cite{horn2012matrix}, the eigenvalues of the matrix $G_n$ lie within its Geshgorin discs, corresponding to its columns (or equivalently to its rows). Subsequently, $G_n$ has eigenvalues with non negative real parts which immediately leads to the fact that it is positive semidefinite \cite{horn2012matrix}. \ill{Condition} (\ref{eq:condition1}) is therefore satisfied. $\square$
\end{proofl}

\begin{prooft}
We will use the dynamics (\ref{eq:model}) and (\ref{eq:network4}) together \il{with
%the passivity conditions in
Assumptions} \ref{assumption:passivedyn2}-\ref{assumption:localminima} to prove Theorem \ref{theorem:converge}. \icl{In particular, it will be shown that the storage functions that follow from the passivity property can be used to construct a Lyapunov function for the network, with stability then deduced using Lasalle's theorem.}

Since the passivity conditions for \il{the} bus dynamics are \illl{considered about the} \il{equilibrium \illl{point},
we
%, the same approach will be adopted for the network system. We therefore
define the deviations from the \illl{corresponding} equilibrium values $\hat{I}_a , \hat{I}_b, \hat{x}, \hat{V}_a, \hat{V}_b$ as $\tilde I_{a,i} = I_{a,i} - \hat{I}_{a,i}, \ \tilde I_{b,i} = I_{b,i} - \hat{I}_{b,i}, \ \tilde x_i =  x_i - \hat{x}_i, \ \tilde V_{a,i} = V_{a,i} - \hat{V}_{a,i}, \ \tilde V_{b,i} = V_{b,i} - \hat{V}_{b,i}$ for the net current injection components, the states and the bus voltage components \il{respectively}.}
%For simplicity, \il{in the remainder of the proof} we drop $\Delta$ and let the variables %$I_{a,i}, \ I_{b,i}, \ x_i, \ V_{a,i}$ and $V_{b,i}$ denote the deviations around the %aforementioned equilibrium.

We now consider the following candidate Lyapunov function for the closed-loop system (\ref{eq:network4}) and (\ref{eq:model}) \il{$\mathcal{V} (x) = \sum_{i=1}^{|\mathcal{N}|} \mathcal{V}_i(x_i)$, where $\mathcal{V}_i(x_i)$ is the storage function
%\footnote{\vsp{The storage function $\mathcal{V}_i(x_i)$ can be used as a valid candidate %Lyapunov function since its time derivative is not increasing with time. Thus, it can be %employed in Lyapunov stability theorem which in combination with La Salle's Invariance %Principle can prove the asymptotic stability of the system.}}
of the \ill{bus dynamics} %(\ref{eq:model})
with input $u_i=[-\tilde I_{a,i}, -\tilde I_{b,i}]$ and output $y_i=[\tilde V_{a,i}, \tilde V_{b,i}]$.}
%\begin{equation*}
%\mathcal{V} (V_a, V_b, x) = \sum_{i=1}^{|\mathcal{N}|} \mathcal{V}_i^B
%\end{equation*}
%which we aim to use later in Lasalle's theorem.
Considering the passivity conditions described in Assumption \ref{assumption:passivedyn2}, we calculate the derivative of the above Lyapunov function \il{with respect to time}. \il{In particular, we get
\begin{equation} \label{eq:derivative1}
\small{
\begin{split}
\dot{\mathcal{V}} = & \sum_{i=1}^{|\mathcal{N}|} \il{\dot{\mathcal{V}}_i} \leq \sum_{i=1}^{|\mathcal{N}|} \Big( [-\tilde I_{a,i} \ -\tilde I_{b,i}] \ \begin{bmatrix} \tilde V_{a,i} \\ \tilde V_{b,i} \end{bmatrix} - \phi_i ( -\tilde I_{a,i}, -\tilde I_{b,i})  \Big) \\
=& - [\tilde V_a^\textrm{T} \ \tilde V_b^\textrm{T}]  H_{2n} \begin{bmatrix} \tilde V_a \\ \tilde V_b \end{bmatrix} - \sum_{i=1}^{|\mathcal{N}|} \phi_i ( -\tilde I_{a,i}, -\tilde I_{b,i})
\end{split}}
\end{equation}
whenever $(\tilde V_{a,i},\tilde V_{b,i}) \in U_i$ and $\tilde x_i \in X_i$ for all $i \in \mathcal{N}$.} Since the matrix $H_{2n}$ and the \il{scalar valued} functions $\phi_i$ are positive semidefinite and positive definite respectively, the inequality (\ref{eq:derivative1}) becomes $\dot{\mathcal{V}} \leq 0$.
%\begin{equation} \label{eq:derivative2}
%\dot{\mathcal{V}} \leq 0.
%\end{equation}
We then \il{make} use of LaSalle's theorem to prove the asymptotic convergence of the system's trajectories to the equilibrium point. According to Assumption \ref{assumption:localminima}, the candidate Lyapunov function $\mathcal{V}$  has \ill{a strict local minimum} at the equilibrium \il{$\hat x$.}
%, which we define as $\hat{\Psi}:=(\hat{I}_a, \hat{I}_b, \hat{x}, \hat{V}_a, \hat{V}_b)$.
\illl{Therefore for a sufficiently small $\epsilon>0$ there exists a compact positively invariant set $\Xi:=\{ x : {\mathcal{V}(x)-\mathcal{V}(\hat x)} \leq \epsilon, \ \hat x \in \Xi, \ \Xi \ \text{connected}\}$ that lies in the neighborhoods stated in the assumptions.}
%By selecting $\epsilon > 0$ sufficiently small \illl{we then choose a  compact positively invariant set $\Xi=\{ x : {\mathcal{V}} \leq \epsilon \}$ such that $\hat x \in \Xi$
%%so \il{that} the states
%and $\Xi$ lies in the neighborhoods stated in the assumptions.}
%we determine the \il{compact positively invariant set $\Xi:=\{ x : {\mathcal{V}} \leq %\epsilon \}$ for the \il{interconnected} system} (\ref{eq:model}) and (\ref{eq:network4}) %that is guaranteed to contain \il{$\hat x$.}
%%$\hat{\Psi}$.
Lasalle's Invariance Principle can now be applied with the function $\mathcal{V}$ on the compact \il{positively} invariant set $\Xi$. This guarantees that all solutions of the \il{interconnected} system (\ref{eq:model}) and (\ref{eq:network4}) with initial conditions $ x(0) \in \Xi$ converge to the largest invariant set within %$\Upsilon:=\Xi \cap \{ \big( I_a, I_b,x, V_a, V_b \big):\dot{\mathcal{V}} = 0\}$
\il{$\Upsilon:=\Xi \cap \{x:\dot{\mathcal{V}} = 0\}$}.
\il{From the positive definiteness of function $\phi$ we have that $\dot{\mathcal{V}} = 0$ implies that $ \tilde I_{a,i} =  \tilde I_{b,i} = 0$, i.e. $I_{a,i}=\hat I_{a,i}$, $I_{b,i}=\hat I_{b,i}$.
Hence from Assumption \ref{assumption:inputoutput} we have that the only \illl{invariant set} in $\Upsilon$ \illl{is} the equilibrium point $x(t)=\hat x$.}
%which under such selection of the positive constant $\epsilon$, contains only the %equilibrium.
\illl{Therefore, for any initial condition $ x(0) \in \Xi$ we have convergence to the equilibrium point, which completes the proof.}
%Finally choosing \il{the set $S$ in the Theorem statement to be} any open neighborhood of %$\hat{x}$ within $\Xi$ completes the proof.
$\square$
\end{prooft}

%\section*{ACKNOWLEDGMENT}

%Acknowledgement goes here.

%%%%%%%%%%%%%%%%%%%%%%%%%%%%%%%%%%%%%%%%%%%%%%%%%%%%%%%%%%%%%%%%%%%%%%%%%%%%%%%%

%%%%% References %%%%%
\balance
\bibliography{report}   %>>>> bibliography data in report.bib
\bibliographystyle{ieeetr}   %>>>> makes bibtex use spiebib.bst

%\begin{IEEEbiography}{Chrysovalantis Spanias}
%received the diploma in Electrical and Computer Engineering from the National Technical University of Athens (NTUA) in 2010. During 2010-13 was employed as a Senior Electrical Engineer in the Construction Industry in Cyprus. Since 2013 he is a PhD student at the Department of Electrical Engineering, Computer Engineering and Informatics, Cyprus University of Technology. His main area of research is stability analysis and control of power systems.
%\end{IEEEbiography}

%\begin{IEEEbiography}{Ioannis Lestas}
%received the B.A. (Starred First) and M.Eng. (Distinction) degrees in Electrical Engineering and Information Sciences and the Ph.D. degree in control engineering from the University of Cambridge (Trinity College) in 2002 and 2007, respectively. His doctoral work was performed as a Gates Scholar and Trinity College Research Scholar. In 2006 he was elected to a Junior Research Fellow of Clare College, University of Cambridge and he was awarded a five year Royal Academy of Engineering research fellowship in 2008. He is also the recipient of an ERC starting grant. He is currently a University Lecturer at Cambridge, Department of Engineering. His research interests include analysis and control of large scale networks with applications in power systems and smart grids, data networks, multiagent systems and biological networks.
%\end{IEEEbiography}

\end{document}